\documentclass[a4paper]{amsart}

\usepackage{amssymb,amsmath,latexsym,amsthm}

\newtheorem{theorem}{Theorem}[section]
\newtheorem{lemma}{Lemma}[section]
\newtheorem{proposition}{Proposition}[section]
\newtheorem{cor}{Corollary}[section]
\theoremstyle{definition}

\newcommand{\beq}[1]{\begin{equation}\label{#1}}
\newcommand{\eeq}{\end{equation}}
\newcommand{\bep}{\begin{proof}}
\newcommand{\eep}{\end{proof}}
\newcommand{\bec}[1]{\begin{cor}\label{#1}}
\newcommand{\eec}{\end{cor}}
\newcommand{\bepr}[1]{\begin{proposition}\label{#1}}
\newcommand{\eepr}{\end{proposition}}
\newcommand{\bel}[1]{\begin{lemma}\label{#1}}
\newcommand{\eel}{\end{lemma}}
\newcommand{\bet}[1]{\begin{theorem}\label{#1}}
\newcommand{\eet}{\end{theorem}}

\def\QQ{\mathbb Q}
\def\ZZ{\mathbb Z}

\def\NN{\mathbb N}

\newcommand{\card}{\mathrm{Card}\kern.7pt}

\let\eps=\varepsilon

\def\hh{{\rm h}\kern.4pt}
\newcommand{\veps}{\ensuremath{\overline{\eps}}}

\newcommand{\vga}{\ensuremath{\overline{\gamma}}}

\begin{document}

\title{There is no Diophantine $D(-1)$--quadruple}

\author[N. C. Bonciocat]{Nicolae Ciprian Bonciocat}
\address{Simion Stoilow Institute of Mathematics of the 
Romanian Academy, Research unit nr. 7,
P.O. Box 1-764, RO-014700 Bucharest, Romania}
\email{Nicolae.Bonciocat@imar.ro}

\author[M. Cipu]{Mihai Cipu}
\address{Simion Stoilow Institute of Mathematics of the 
Romanian Academy, Research unit nr. 7,
P.O. Box 1-764, RO-014700 Bucharest, Romania}
\email{Mihai.Cipu@imar.ro}

\author[M. Mignotte]{Maurice Mignotte}
\address{Universit\'e de Strasbourg, U. F. R. de Math\'ematiques,
 67084 Strasbourg, France}
\email{Maurice.Mignotte@math.unistra.fr}

\thanks{This joint work has been initiated as a  project  
LEA Franco-Roumain Math-Mode and completed in the framework 
of a project GDRI ECO-Math. The second author is grateful for 
financial support allowing to attend ANTRA 2017 Conference
organized by RIMS Kyoto and XVI Conference on Representation Theory,
Dubrovnik, 2019, where  preliminary results were discussed.
}
\date{\today}%

\begin{abstract}
A set of  positive integers with the property that the 
product of any two of them is the successor of a perfect 
square is called Diophantine $D(-1)$--set. Such objects are 
usually studied via a system of generalized Pell equations 
naturally attached to the set under scrutiny. In this 
paper, an innovative technique is introduced in the study 
of Diophantine $D(-1)$--quadruples. The main novelty is the 
uncovering of a quadratic equation relating various parameters 
describing a hypothetical $D(-1)$--quadruple with integer 
entries. In combination with extensive computations, this 
idea leads to the confirmation of the conjecture according
to which there is no  Diophantine $D(-1)$--quadruple.
\end{abstract}

\maketitle

\section{The strategy}\label{sec1}

In the third century, Diophantus of Alexandria found four 
positive rationals such that the product of any two of them
increased by unity is a square, see, for 
instance,~\cite{dic,diog,dior,hea}. 
Fermat found the quadruple  consisting of positive integers 
$1$, $3$, $8$, $120$ with the same property. As Euler remarked, Fermat's
set can be enlarged by inserting $\frac{777480}{8288641}$
without losing the defining property. It was in 1969 that
Baker and Davenport~\cite{bad} proved that there is no quintuple of
positive integers containing Fermat's set and still having the
property of interest. On this occasion the authors introduced 
an important tool, nowadays referred to as Baker--Davenport 
lemma, for the effective resolution of Diophantine equations.

Diophantus also studied a problem that turned out to be closely
related to that mentioned before. Namely, he asked for numbers
such that the product of any two of them increased by the sum
of these two is a square. Since $ab+a+b=(a+1)(b+1)-1$, the
question boils down to finding sets with the property that
the product of any two of its elements is one more a square.
The essence of both problems is captured by the next definition.

Let $m\ge 2$ and $n$ be integers. A set of $m$ positive integers 
is called Diophantine $D(n)$--$m$--set if the product of any two 
distinct elements increased by $n$ is a perfect square. In this terminology, 
Fermat's example is a $D(1)$--quadruple, and the set $\{ 4, 9, 28 \}$ 
presented by Diophantus himself as an answer to the second problem 
gives rise to the $D(-1)$--triple $\{5,10,29\}$. A more general 
notion is obtained by considering elements of any  commutative 
ring  instead of positive integers. However, many difficult, 
interesting problems  already occur in the setting fixed by the 
given definition. In the rest of the paper we shall refer only
to this definition, even when we omit the adjective ``Diophantine''.

It is worth mentioning that the objects produced by this definition 
with  $n=0$ are not particularly interesting --- for each positive
$m$ there exist infinitely many $D(0)-m$--sets and even infinite 
$D(0)$--sets. Therefore, when speaking of $D(n)-m$--sets we shall 
always assume $n$ is nonzero.

A natural question is how large a $D(n)$--set can be. 
 It is known~\cite{duca} that for $\vert n\vert <400$ 
one has $m<31$, and for other $n$ the cardinality of any 
$D(n)$--$m$--tuple is at most $16 \log \vert n\vert$. Better 
bounds are known for particular values of $n$. As noticed
in several papers (among which~\cite{br,gusi,mora}), for $n\equiv
2 \pmod 4$ there is no $D(n)$--quadruple. On the opposite side,
in~\cite{du4} it is showed that if  $n \not \equiv 2 \pmod 4$ 
and $n\not \in S:=\{ -4,-3,-1,3,5,8,12,20\},$ then there exists
at least one $D(n)$--quadruple. In the same paper Dujella 
expressed his confidence that this is all one can hope for.

\medskip
\noindent
\textbf{Conjecture.} There exists no  $D(n)$--quadruple for 
$n\in \{ -4,-3,-1,3,5,8,12,20\}$.

\medskip

According to a remarkable result of Dujella and Fuchs~\cite{dufu},
in any $D(-1)$--quadruple $(a,b,c,d)$ with $a<b<c<$ one has $a=1$.
This readily implies the nonexistence of  $D(-1)$--quintuples. 
The same authors together with Filipin proved in~\cite{duff} that
there are at most finitely many $D(-1)$--quadruples. 
The present authors obtained in~\cite{bcm} the bound $10^{71}$
for the number of $D(-1)$--quadruples, thus improving on the
previous bound $10^{356}$ found in~\cite{fifu}.  Better estimates
have been given lately: in~\cite{eff} one finds the upper bound
$5\cdot 10^{60}$, successively strenghtened to $3.01\cdot 10^{60}$ 
in~\cite{tim}, to $4.7\cdot 10^{58}$ in~\cite{lap1}. The best 
bound we are aware of is $3.677\cdot 10^{58}$ found in~\cite{lap2}.

A basic technique in the study of $D(n)$--sets exploits a connection
with systems of generalized Pell equations. We explain the main
ideas of this approach in the framework of $D(-1)$--quadruples.

Suppose $(1,b,c,d)$  is a  $D(-1)$--quadruple with $1<b<c<d$.
Then there are positive integers $r$, $s$, $t$, $x$, $y$, $z$
satisfying
\beq{ecbc}
b-1=r^2, \ c-1=s^2, \ bc-1=t^2, 
\eeq
\beq{ecd}
d-1=x^2, \ bd-1=y^2, \  cd-1=z^2. 
\eeq
Eliminating $d$ in Eq.~\eqref{ecd}, one obtains a system of 
three generalized Pell equations 
\beq{ecpc}
z^2-cx^2=c-1, 
\eeq 
\beq{ecpbc}
bz^2-cy^2=c-b,
\eeq
\beq{ecpb}
y^2-bx^2=b-1.
\eeq
By Theorem~1.2 in~\cite{bcm}, we may assume $c<2.5 \, b^6$. Then, 
according to~\cite[Lemmata~1 and~5]{duff}, the positive integer 
solutions of each of the above Pell equations are respectively 
given by
\[
 z+x\sqrt{c} =s(s+\sqrt{c})^{2m}, \quad m\ge 0,
\]
\[
z\sqrt{b}+y\sqrt{c}=(s\sqrt{b}+\rho r\sqrt{c}) 
(t+\sqrt{bc})^{2n},\quad  n\ge 0,
\]
\[
y+x\sqrt{b}=r(r+\sqrt{b})^{2l}, \quad l\ge 0,
\]
for   fixed $\rho \in \{-1,1 \}$.
Therefore, the triples  $(x,y,z)$ of positive integers that
simultaneously satisfy Eqs.~\eqref{ecpc}--\eqref{ecpbc} are such that
\beq{ec12mn}
 z=v_m=w_n,
 \eeq
where the integer sequences $(v_p)_{p\ge 0}$, $(w_p)_{p\ge 0}$ 
are given by  explicit formul\ae \ 
\beq{eqv}
 v_{p}=\frac{s}{2} \left( (s+\sqrt{c})^{2p}+
(s-\sqrt{c})^{2p}\right)
\eeq
and respectively
\beq{eqw}
w_p=\frac{s\sqrt{b} +\rho r\sqrt{c}}{2\sqrt{b}}(t+\sqrt{bc})^{2p}
+\frac{s\sqrt{b} -\rho r\sqrt{c}}{2\sqrt{b}}(t-\sqrt{bc})^{2p}.
\eeq

These formul\ae \ give rise in the  usual way to linear forms 
in the logarithms of three algebraic numbers, for which upper
bounds are obtained directly, while lower bounds are given by
a general theorem of Matveev~\cite{mat}. Comparison of these
bounds results in inequalities for indices $m$ and $n$ in terms
of elementary functions in $b$ and $c$. In order to get reverse
inequalities, Dujella and Peth\H o introduced in~\cite{dup}
the congruence method. Their idea is to consider the recurrent
sequences modulo $8c^2$ and prove that suitable hypotheses entail
that these congruences are actually equalities. The best result
obtained by this approach is due to Dujella, Filipin and
Fuchs~\cite{duff}.

\bet{tedff} 
Let $(1,b,c,d)$ with $1<b<c<d$ be a  $D(-1)$--quadruple. 
Then $b>100$ and $ c< \min \{ 11\, b^6, 10^{491} \}$. 
More precisely:
\begin{enumerate} 
\item [a)] If $b^3\le c < 11 \, b^6$, then $c< 10^{238}$.
\item [b)] If $b^{1.1} \le c< b^{3}$, then $c< 10^{491}$.
\item [c)] If $3b \le c< b^{1.1}$, then $c< 10^{94}$.
\item [d)] If $b < c< 3b$, then $c< 10^{74}$.
\end{enumerate}
\eet

A variant of the congruence method has been introduced 
in~\cite{bcm}. The new idea is to interpret an equivalence 
$L \equiv R \pmod c$ as an equality $L-R=jc$ for a suitable 
integer $j$. Instead of striving to get $j=0$, as did  the 
predecessors, all possibilities for the sign of $j$ have 
been analysed. As a result of this study,
inequalities of the form $n>f(b,c)^{\alpha (j)}$ have been
established. Combined with another new idea, called 
smoothification in~\cite{bcm}, and large-scale computations,
always performed with the help of the package 
PARI/GP~\cite{pari2}, this yields much better results.

\bet{tenoi}
Let $(1,b,c,d)$ with $1<b<c<d$ be a  $D(-1)$--quadruple. 
Then $b>1.024\cdot 10^{13}$ and $ \max \{10^{14}b,  b^{1.16} \}
 < c< \min \{ 2.5\, b^6, 10^{148} \}$. 
More precisely:
\begin{enumerate}
\item [i)] If $b^5\le c < 2.5 \, b^6$, then $c< 10^{100}$.
\item [ii)] If $b^4\le c <b^5$,  then $c< 10^{82}$.
\item [iii)] If $b^{3.5}\le c< b^4$, then $c< 10^{66}$.
\item [iv)] If $b^3\le c< b^{3.5}$, then $c< 10^{57}$.
\item [v)] If $b^{2}\le c< b^3$, then $c< 10^{111}$. 
\item [vi)] If $b^{1.5}\le c< b^2$, then $c< 10^{109}$.
\item [vii)] If $b^{1.4} \le c< b^{1.5}$, then $c< 10^{128}$.
\item [viii)] If  $b^{1.3} \le c< b^{1.4}$, then $c< 10^{148}$.
\item [ix)] If    $b^{1.2} \le c< b^{1.3}$, then $c< 10^{133}$.
\item [x)] If $b^{1.16} \le c< b^{1.2}$, then $c< 10^{107}$.  
\end{enumerate}
\eet

More recently, Filipin and Fujita obtained an even better
relative bound for the third element of a hypothetical
 $D(-1)$--quadruple. In~\cite{ff} they proved the remarkable
result quoted below. The proof is based on an improved variant
of Rickert's theorem~\cite{ric}.

\bet{tefifu}
Any  $D(-1)$--quadruple $(1,b,c,d)$ with $1<b<c<d$ satisfies
$c< 9.6 \, b^4$.
\eet

\begin{figure}[t]
\begin{picture}(344,160)\setlength{\unitlength}{1pt}
\put(5,0){\vector(1,0){344}} %
\multiput(72,-1)(67,0){5}{\line(0,1){2}}   %
\put(64,-11){\makebox(20,10){{\small 100}}}    %
\put(131,-11){\makebox(20,10){{\small 200}}}  
\put(198,-11){\makebox(20,10){{\small 300}}}  
\put(265,-11){\makebox(20,10){{\small 400}}}  
\put(312,-11){\makebox(20,10){{\small $\log _{10} c$}}}  
\put(5,0){\vector(0,1){165}} %
\multiput(4,0)(0,30){6}{\line(1,0){2}}   %
\put(-5,-5){\makebox(10,10){{\small 1}}}    %
\put(-5,25){\makebox(10,10){{\small 2}}} 
\put(-5,55){\makebox(10,10){{\small 3}}} 
\put(-5,85){\makebox(10,10){{\small 4}}} 
\put(-5,115){\makebox(10,10){{\small 5}}} 
\put(-5,145){\makebox(10,10){{\small 6}}} 
\put(-20,130){\makebox(20,10){{\small $\log _b c$}}} 
\put(5,150){\line(1,0){160}}       %
\put(72,150){\line(0,-1){30}}  
\put(60,120){\line(1,0){12}}
\put(60,120){\line(0,-1){30}} 
\put(49,90){\line(1,0){11}}
 \put(49,90){\line(0,-1){15}}
\put(42,75){\line(1,0){7}}
\put(42,75){\line(0,-1){15}}
\put(42,60){\line(1,0){39}}
\put(81,60){\line(0,-1){30}}
\put(77,30){\line(1,0){4}}
\put(77,30){\line(0,-1){15}}
\put(77,15){\line(1,0){12}}
\put(89,15){\line(0,-1){3}}
\put(89,12){\line(1,0){13}}
\put(102,12){\line(0,-1){3}}
\put(92,9){\line(1,0){10}}
\put(92,9){\line(0,-1){3}}
\put(75,6){\line(1,0){17}}
\put(75,6){\line(0,-1){1}}
\put(5,5){\line(1,0){70}}       %
\put(165,150){\line(0,-1){90}}  %
\put(165,60){\line(1,0){170}}
\put(334,60){\line(0,-1){57}}
\put(68,3){\line(1,0){266}}
 \put(68,3){\line(0,-1){3}}     %
\thicklines 
\put(5,90){\line(1,0){336}}     %
\end{picture}
\caption{Absolute bounds for $c$ given by Theorems~\ref{tedff} 
to~\ref{tefifu}}\label{fig1}
\end{figure}
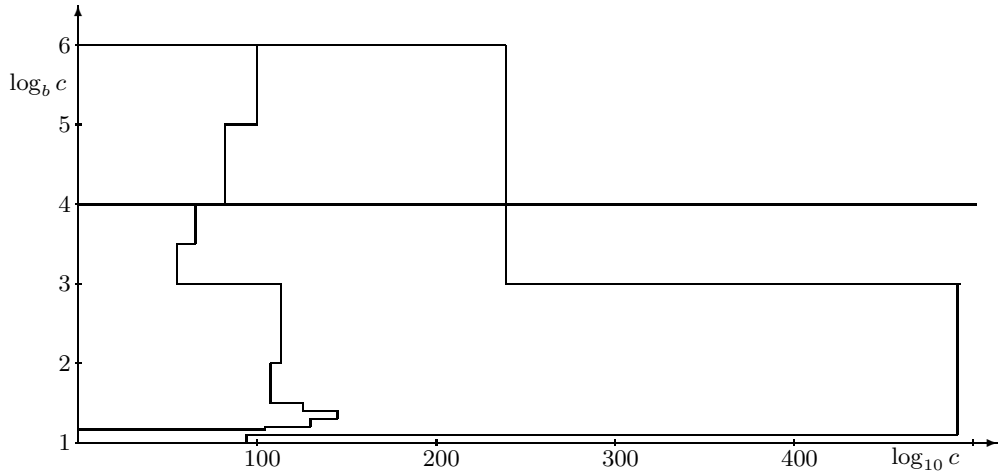

\begin{figure}[b]
\begin{picture}(330,150)\setlength{\unitlength}{1pt}
\put(-5,0){\line(1,0){330}}
\put(234,0){\circle*{2}}           
\put(234,0){\line(2,1){46}}
\put(280,23){\circle*{2}}          
\put(280,23){\line(5,3){32}}
\put(312,42){\circle*{2}}          
\put(312,42){\line(3,4){17}}
\put(329,65){\circle*{2}}          
\put(329,65){\line(-6,1){39}}
\put(290,72){\circle*{2}}          
\put(290,72){\line(-3,2){62}}
\put(228,113){\circle*{2}}         
\put(228,113){\line(-4,1){101}}
\put(127,138){\circle*{2}}         
\put(127,138){\line(-1,-4){18}}
\put(109,66){\circle*{2}}          
\put(109,66){\line(-6,-1){35}}
\put(74,60){\circle*{2}}          
\put(74,60){\line(-2,-1){52}}
\put(22,34){\circle*{2}}           
\put(22,34){\line(-3,-5){20}}
\put(2,0){\circle*{2}}           
\thicklines
\put(146,0){\circle{3}}            
\put(146,0){\line(-5,4){109}}
\put(50,85){\makebox{{\small $c=9.6 \, b^4$}}}
\put(8,-10){\makebox{{\small $c=b^6$}}}
\put(274,-10){\makebox{{\small $c=b^{1.16}$}}}
\end{picture}
\caption{A different view on the absolute bounds for $c$ given by 
Theorem~\ref{tenoi}}\label{fig2}
\end{figure}
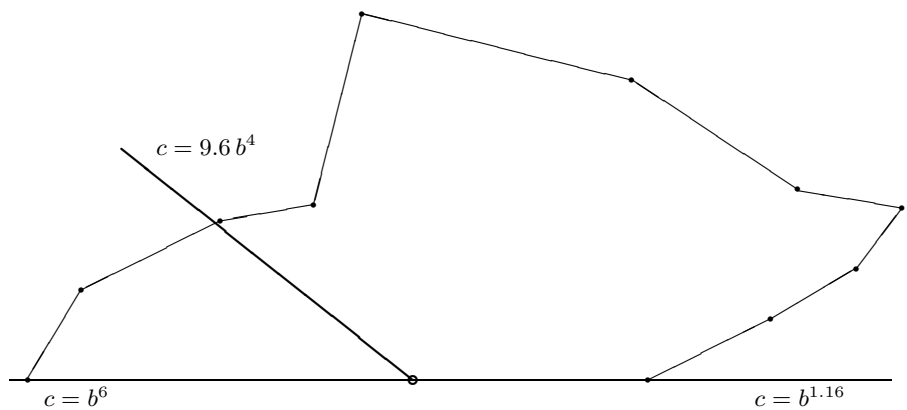

Figures~\ref{fig1} and~\ref{fig2} give a graphical representation
of these results.
The outer (inner) polygon in Figure~\ref{fig1} represents the
region where $c$ is confined according to Theorem~\ref{tedff}
(Theorem~\ref{tenoi}). From Theorem~\ref{tefifu} we see that
in the search of $D(-1)$--quadruples we can restrict ourselves
to the part of polygon sitting in the half-plane $c< 9.6 \, b^4$.
Figure~\ref{fig2} contains an approximate illustration of
Theorem~\ref{tenoi} in polar coordinates. It is seen that $c$
is to be found in a region whose shape looks like a nonstandard 
fan. We interpret the presence of  inlets as a strong hint 
that actually there are no $D(-1)$--quadruples whose third entry 
is located in the nonconvex blades. Eliminating these blades has
the effect of ``partially closing the fan''. The aim of this paper 
is to ``completely close the fan''. In this respect we will prove 
the following result.

\bet{tede}
There is no Diophantine $D(-1)$--quadruple.
\eet

Since, by~\cite[Remark~3]{du4}, all elements of a $D(-4)$--quadruple 
are even, from Theorem~\ref{tede} we get for free another result
that provides partial confirmation of Dujella's conjecture.

\bet{tede4}
There is no Diophantine $D(-4)$--quadruple.
\eet

Our strategy is based  on the following interpretation of
Eq.~\eqref{ecbc}: the initial triple $(1,b,c)$ of a hypothetical
$D(-1)$--quadruple with $1<b<c<d$ is associated to a member of a
two-parameter family of integers $b=r^2+1$, $c=s^2+1$ and $t$ 
witnesses the fact that we deal with a $D(-1)$--triple. We would 
like to handle all these parameters simultaneously. This goal is 
achieved by considering the integer
\[
 f=t-rs,
\]
which is easily seen to be positive. It already appeared in several
proofs available in literature, see, for instance,~\cite{dufu,fuaa,furo,het}.
Up to our work it always had a secondary role, sitting in
background. Focusing on $f$ turned out to open new prospects for 
the study of $D(-1)$--quadruples.

The developments below seem to have been overlooked in the
literature.
Squaring  $f+rs=t$, one gets $f^2+2frs+r^2s^2=(r^2+1)(s^2+1)-1$,
whence
\beq{ecrs} 
r^2+s^2=2frs+f^2.
\eeq

Our approach hinges on the study of solutions in positive 
integers to the master equation~\eqref{ecrs}  in its various 
disguises. This study is much easier than the examination of 
solutions to the system of generalized Pell equations~\eqref{ecpc} 
to \eqref{ecpb}. As will be seen in Section~\ref{sec2}, rather
strong results are obtained by elementary proofs relying
on properties of solutions to equation~\eqref{ecrs}. 
Besides the emphasis on $f$ already mentioned, our treatment
introduces here yet another variant of the congruence method
by considering  modulus $8f^6$ instead of the ``classical'' $8c^2$.
Also, we revisit published results and use them in a novel way. 
The study is straightforward if $\gcd (r,s)=f$.  
In the general case, we complement it with
considerations along the lines described in the next
paragraphs.

Let us denote by $J$ the set of pairs of integers $(r,s)$ such 
that there exists a $D(-1)$--quadruple $(1,b,c,d)$ with
$1<b<c<d$ and $b=r^2+1$, $c=s^2+1$. For $(r,s)\in J$ we put  
$s=r^\theta$ and define
\[
\theta^-=\inf_{(r,s)\in J} \theta, \quad \theta^+=\sup_{(r,s)\in J}
\theta.
\]
These numbers measure the size of $c$ with respect to $b$.
``Closing the fan'' means showing that for a putative 
$D(-1)$--quadruple one has $\theta^+ < \theta^-$.

The main result of~\cite{bcm} (quoted in Theorem~\ref{tenoi}) 
gives in particular $c\le 3 \, b^6$. Using this upper bound, 
a computer verification described in~\cite[Section~2]{bcm}  
led to the conclusion that $r>32\times 10^5$, so that
$b>10^{13}$. Hence
\[
  \theta^+ \le 6{.}1.
\]
The approach followed in~\cite{bcm} used, among other things,
the inequalities $c>3.999 \, f^2 b$ and $f> 10^7$ in order
to increase the lower bound $\theta^-$ from $1$ to $1.16$.
Pursuing this idea requires to examine much higher values
of $f$, a process that becomes prohibitively time-consuming.
To give the reader a feeling of the difficulties encountered,
we mention that computations that allow exclusion of values 
$f \le 10^7$ needed about two weeks (measured by wall-clock)
on a personal computer; to reach the level
$f > 10^8$, our program ran on a network of up to 6 computers
for about three months; further computations were performed
during other six months on as much as 30 computers.
Therefore, another course has been chosen: instead of shortening
the interval $[\theta^-,  \theta^+]$, we looked for methods of
splitting it in such a way that parallel processing is possible.
The alternative approach was devised after it was observed that
there is no $D(-1)$--quadruple with $b^2 \le c \le 16\, b^3$
even when $\gcd (r,s) \ne f$  and requires to identify suitable 
functions of parameters already introduced in the formulation of 
the problem.  The breakthrough was realized after
making the choice $F:=s-2rf$. A short study revealed that 
it is very helpful in separating values $c < 4b^2$ from values
$c > 4b^2$. It was a matter of days to reach the conclusion
that around $c =4b^2$ there is a large gap (a comparatively  
long interval in which there is no third entry of a
$D(-1)$--quadruple). 

Completion of the proof requires new explicit computations. Besides 
those needed for use of results providing bounds for linear forms
in logarithms, a considerable amount of them was devoted to solve
many quadratic Diophantine equations and then to apply the reduction
procedure based on Baker--Davenport lemma. Further explanations and 
full details are given in Section~\ref{sec4}.

A successful implementation of the strategy just sketched
requires to pay attention to several aspects which will be
described in Section~\ref{secAleks}. Here we mention only one point.
It is clear that the smaller the upper bounds on $b$ and $c$
are, the faster the subsequent computations depending on them
are. To this end,  Matveev's general theorem~\cite{mat} was first 
replaced by a strengthening of it due to Aleksentsev~\cite{alex} 
and next by an older result of  Matveev~\cite{Mat98}. 
This course of action is determined by our experience, according 
to which a giant step is better replaced by succesive small steps.

\section{A sufficient condition for the nonexistence of
$D(-1)$--quadruples}\label{sec2}

The aim of this section is to revisit results of our previous
work~\cite{bcm} in the light of the new guiding strategy. As it
turns out, a lot of information already available can be exploited 
in a novel manner, producing unexpected results and suggesting 
further developments. 

The starting point of our study of solutions in  positive integers 
to the equation 
\[ 
r^2+s^2=2frs+f^2,
\]
is the observation that $\gcd (r,s)$ is a divisor of $f$. As 
extremal elements / circumstances generally are very interesting,
we consider solutions to Eq.~\eqref{ecrs} such that  
\beq{ecmild}
f=\gcd (r,s).
\eeq
Then  one has
\[
r=fu, \quad s=fv,  
\]
for some positive integers  $u$, $v$ satisfying
\[
(v-fu)^2-eu^2=1,
\]
with $e=f^2-1$.

Let $\gamma=f+\sqrt{e}$ be the fundamental solution to the Pell
equation $V^2-eU^2=1$ and $\vga=f-\sqrt{e}$ its  algebraic 
conjugate. According to Lemma~3.5 from~\cite{bcm}, 
$v>3.999^{1/2}f u$, so that all positive solutions to 
Eq.~\eqref{ecrs}  have the form
\beq{ec1}
r=\frac{f(\gamma ^k-\vga ^k)}{\gamma -\vga},
\quad s=\frac{f(\gamma ^{k+1}-\vga ^{k+1})}{\gamma -\vga},
\quad k\in \NN,
\eeq
whence
\beq{ec2}
b=\frac{f^2(\gamma ^k+\vga ^k)^2-4}{(\gamma -\vga)^2},
\quad 
c=\frac{f^2(\gamma ^{k+1}+\vga ^{k+1})^2-4}{(\gamma -\vga)^2}.
\eeq

Notice that the main result of He--Togb\'e~\cite{het} assures
$f\ge 2$, a piece of information we shall repeatedly use 
without explicitly mentioning it. Later on, a more stringent
restriction, checked computationally, will be preferred.

With this notation fixed, we proceed to examine the properties
of solutions to Eq.~\eqref{ecrs} under the 
condition~\eqref{ecmild}.

\bepr{pr21p}
One has $c<b \, \gamma ^2$. Moreover, for $k\ge 2$ it holds
\[
\gamma ^2-\frac{1}{2}<\frac{c}{b}. 
\]
\eepr \bep
The numerator of $b\gamma ^2-c$ is $(\gamma ^2-1)[f^2(2+
\vga^{2k}+  \vga^{2k+2})-4]$, which is manifestly positive. 
The numerator of $2c-(2\gamma ^2-1)b$  is found to be 
\[
 4(2\gamma ^2-3)+f^2[ \gamma ^{2k}-2(2\gamma ^2-3)
+ \vga ^{2k}(1-2\gamma ^2+2\vga ^2)].
\]
Since $k\ge 2$ and  $0<\vga <1$, it is sufficient to prove that  
\[
\gamma ^4 > 2(2\gamma ^2-3)+2\gamma ^2-1. 
\]
This inequality is obvious  on noticing the identity
$\gamma ^4=4(f^2-1)\gamma ^2+2\gamma ^2-1$.
\eep

We can bound from below $b$  by a power of $\gamma$.

\bepr{pr22p}
If $(1,b,c,d)$ is a $D(-1)$--quadruple with $10^{13} <b<c<d$
and $b$, $c$ given by formula~\eqref{ec2}, then
\[
 \gamma ^{2k-1} <b.
\]
\eepr \bep
The desired inequality is equivalent to $f^2(\gamma ^{2k}+2+ 
\vga^{2k})-4 > 4(f^2-1) \gamma ^{2k-1}$, which follows 
from $ f^2\gamma > 4(f^2-1)$.
\eep

We are now in a position to show that the third entry in a
$D(-1)$--quadruple restricted as in~\eqref{ecmild} is much 
closer to the second one than was previously known.

\bepr{pr23p}
Any $D(-1)$--quadruple $(1,b,c,d)$ with $10^{13} <b<c<d$  and 
$b$, $c$ given  by~\eqref{ec2} satisfies $c<b^3$.
\eepr \bep
Assuming $b^3\le c$, we deduce with the help of the previous 
results 
\[
\gamma ^{4k-2}<b^2\le \frac{c}{b}<\gamma^2, 
\]
that is, $k<1$. Thus  $k=0$, whence $r=0$ and $b=1$, which 
is not possible.
\eep

Propositions~\ref{pr21p} and~\ref{pr22p} have other important
consequences drawn from information made available by our previous
work.
For the sake of convenience, we recall an experimental
result obtained after two weeks of computer calculations
for the needs of Lemma 3.5 from~\cite{bcm}. Its proof is based
on the well-known structure of solutions to a Pellian equation
of the type
\beq{ecf1}
W^2-D U^2=N, \quad D >0 \quad  \mathrm{nonsquare \ and} \quad  N\ne 0.
\eeq
The most familiar reference is Nagell's book~\cite{nag} in its 
various editions but the results have been published already in 
the 19th century by Chebyshev~\cite{ceb}.  More details
are available in the proof of~\cite[Lemma~2.9]{bcm}.

\bepr{prexp} There are no  $D(-1)$--quadruples  $(1,b,c,d)$ 
with the corresponding $f$ less than or equal to $10^7$.
\eepr 

This result is used below in conjunction with the fact
that for any hypothetical $D(-1)$--quadruple one has $n<10^{19}$
(see Table~1 from~\cite{bcm}).

After these preparations, we proceed with the study of positive 
solutions to Eq.~\eqref{ecrs} given by~\eqref{ec1}.

\bepr{pr34} There are no  $D(-1)$--quadruples  $(1,b,c,d)$ 
with $10^{13} <b<c<d$,  $b^2\le c<b^3$   and  $b$, $c$ 
given  by~\eqref{ec2}. 
\eepr \bep
Suppose, by way of contradiction, that the thesis is false.
From 
\[
\gamma ^{2k-1}<b\le \frac{c}{b}<\gamma^2
\]
we get $k<2$. Since $k\ge 1$, we conclude that $k=1$.

 Eliminating $d$ in Eq.~\eqref{ecd} yields the system of 
generalized Pell equations~\eqref{ecpc}--\eqref{ecpb}.
It is well known that $z$ appears in two second-order 
linearly recurrent sequences. Thus (see, e.g.,~\cite{dufu} 
or~\cite{duff}) $z=v_m=w_n$, with $m$, $n$ positive integers 
of the same parity, and
\beq{relv}
v_0=s, \quad  v_1=(2c-1)s, \quad v_{m+2}=(4c-2) v_{m+1} 
-v_{m},
\eeq
\beq{relw}
w_0=s, \quad  w_1=(2bc-1)s +2\rho rtc, \quad w_{n+2}=
(4bc-2) w_{n+1} -w_{n},
\eeq
where $\rho =\pm 1$. Since $k=1$, one has
$r=f$, $b=f^2+1$, $s=2f^2$, $c=4f^4+1$, $t=2f^3+f$, and
therefore
\[
v_{m+2}=2(8f^4+1) v_{m+1} -v_{m},
\] 
\[
w_{n+2}=2(8f^6+8f^4+2f^2+1) w_{n+1} -w_{n}.
\] 
Taken modulo $8f^6$, these recurrent relations readily give
\[
v_m\equiv 2f^2 \pmod{8f^6} 
\]
and  
\[ 
w_{n}\equiv  2(\rho n+1) f^2 + \bigl(\frac{4n^3+8n}{3}\rho +4n^2
\bigr) f^4 \pmod{8f^6} .
\] 
Together with $v_m=w_n$, this implies 
\beq{ec6}
\rho n  + \left(\frac{2n^3+4n}{3}\rho +2n^2 \right) f^2 
\equiv 0 \pmod{4f^4}.
\eeq
Note that $n$ must be even and use this information to deduce
$n\equiv 0 \pmod{4f^2}$,
 so that $n=4f^2u$ for some positive integer $u$. Replace $n$
by $4f^2u$  in~\eqref{ec6} to get $u\equiv 0 \pmod{f^2}$, and
therefore $n\ge 4f^4$.

Using this inequality  and Proposition~\ref{prexp} one gets 
$n> 10^{28}$, in contradiction with the fact $n<10^{19}$ proved
in~\cite[Proposition~4.3]{bcm}.
\eep

\bepr{pr35} There are no  $D(-1)$--quadruples $(1,b,c,d)$ 
with $10^{13} <b<c<d$,  $b^{1.5}\le c<b^2$, and  $b$, $c$ 
given  by~\eqref{ec2}.
\eepr  \bep
We reason by reduction to absurd. Suppose that  $(1,b,c,d)$
is a   $D(-1)$--quadruple satisfying  $10^{13} <b<c<d$,  
$b^{1.5}\le c<b^2$, and  $b$, $c$  given  by~\eqref{ec2}.
It is easy to prove the upper bound $2k< 5$  as above. 
Assuming $k=1$, from $ c<b^2$ it results $4f^4+1<(f^2+1)^{2}$, 
whence $f<1$, which is impossible. %

So it is established that $k= 2$. Therefore, one has
$r=2f^2$, $s=4f^3-f$, and $t=8f^5-2f^3+f$. The solutions to  
the system of Pellian equations~\eqref{ecpc}--\eqref{ecpb} 
verify $y=U_n=u_l$, where $n$, $l$ are positive integers and
\[  u_{l+2}=(4b-2)u_{l+1}-u_l,  \quad u_0=r, \quad u_1=(2b-1)r,
 \]
\[
  U_{n+2}=(4bc-2)U_{n+1}-U_n,  \quad U_0=\rho r, \quad
    U_1=(2bc-1)\rho r +2bst.
\]
Considering these recurrence relations modulo $r^3$, one readily gets
\[
  u_l\equiv r \pmod{r^3}.
\]
A short inductive reasoning that takes into account the explicit
formul\ae \ giving $r$, $s$, $t$ in terms of $f$ results in the 
congruence
\[  
 U_n \equiv (n^2r^2+r)\rho +\frac{(10n-n^3)}{3} r^2-n r \pmod{r^3},
\]
so that $u_l=U_n$ implies $r \rho -nr \equiv r \pmod{r^3}$. Therefore,
there exists an integer $\lambda$ such that
\[
n = \rho -1 + \lambda r.
\]
Since $n\ge 7$ by~\cite[Proposition~2.2]{bcm}, it follows that 
$\lambda$ is positive. Introducing this formula for $n$ in the 
congruence for $U_n$, one sees that in fact  one has $\lambda \ge r-1$, 
so that  $n>0.9 r^2 > 3 f^4 >10^{19}$. As explained previously, 
this contradicts~\cite[Proposition~4.3]{bcm}. The contradiction 
is due to the assumption that there exists a  $D(-1)$--quadruple
satisfying all the hypotheses of the present proposition.
\eep

\bepr{pr36} There are no  $D(-1)$--quadruples $(1,b,c,d)$ 
with $10^{13} <b<c<d$,  $b^{1.4}\le c<b^{1.5}$,
and  $b$, $c$ given  by~\eqref{ec2}.
\eepr \bep
As above, we reason by reduction to absurd. Suppose that  
$(1,b,c,d)$ is a $D(-1)$--quadruple satisfying  
$10^{13} <b<c<d$, $b^{1.4}\le c<b^{1.5}$, and  $b$, $c$ given  
by~\eqref{ec2}. As seen in the proof of the previous  result, 
one then has $k\ge 2$. For $k=2$, from $ c<b^{1.5}$ one gets 
\[
 c=16f^6-8f^4+f^2+1<(4f^4+1)^{1.5}<9f^6+1,
\]
that is, $f< 1$, a contradiction. Therefore, we
conclude that $k\ge 3$. This and Proposition~\ref{pr22p}
yield
\[
\frac{c}{b}\ge b^{0.4} > \gamma ^{0.4(2k-1)} \ge \gamma ^2, 
\]
 which contradicts Proposition~\ref{pr21p}.
\eep

\bepr{pr37} There are no  $D(-1)$--quadruples $(1,b,c,d)$ 
with $10^{13} <b<c<d$,  $b^{1.3}\le c<b^{1.4}$,
and  $b$, $c$ given  by~\eqref{ec2}.
\eepr \bep
From the last proof we retain that $k$ is at least $3$,
while from the chain of inequalities
\[
 \gamma ^2> \frac{c}{b}\ge b^{0.3} > \gamma ^{0.3(2k-1)} 
\]
we deduce that $k\le 3$. To obtain a bound for $f$, we
follow the reasoning in the proof of Proposition~\ref{pr34}.

Since $k=3$, one has
\[
v_m\equiv -4f^2+8f^4 \pmod{8f^6},
\]
\[
w_n\equiv -(2\rho n+4)f^2+\bigl(\frac{4n-4n^3}{3}\rho -8n^2
+8\bigr) f^4 \pmod{8f^6},
\]
whence  again it follows $n\ge 4f^4$. As already seen, this 
leads to a contradiction.
\eep

\bepr{pr38} There are no  $D(-1)$--quadruples $(1,b,c,d)$ 
with $10^{13} <b<c<d$,  $b^{1.2}\le c<b^{1.3}$,
and  $b$, $c$ given  by~\eqref{ec2}.
\eepr \bep
We adapt the reasoning used to establish Proposition~\ref{pr35}.
So let $(1,b,c,d)$ be a $D(-1)$--quadruple satisfying  
$10^{13} <b<c<d$,  $b^{1.2}\le c<b^{1.3}$,
and  $b$, $c$ given  by~\eqref{ec2}.

In the proof of Proposition~\ref{pr36} it was shown that
$k\ge 3$.  For $k=3$, from $c<b^{1.3}$ it results
\[
16f^4(2f^2-1)^2+1 < \left(f^2(4f^2-1)^2+1\right)^{1.3}
< 40f^8+1, 
\]
whence $f^2<3$,  a contradiction. Therefore, one has 
$k\ge 4$. From
\[
 \gamma ^2> \frac{c}{b}\ge b^{0.2} > \gamma ^{0.2(2k-1)} 
\]
we deduce that $k\le 5$. Assuming $k=5$, one obtains 
$n\ge 4f^4$ as  in the proof of Proposition~\ref{pr34},
so a contradiction appears in this case. It remains to examine 
what happens when $k=4$, $r=8f^4-4f^2$, $s=16f^5-12f^3+f$,
and $t=128f^9-160f^7+56f^5-4f^3+f$. 

A short study of the sequences $(u_l)_l$, $(U_n)_n$ introduced in 
the proof of Proposition~\ref{pr35} gives 
\[  
U_n \equiv \bigl( (8-8n^2)f^4 -4f^2 \bigr) \rho +
\frac{(4n^3-100n)}{3} f^4 + 2 nf^2 \pmod{8f^6}.
\]
Proceeding as in the proof  of Proposition~\ref{pr35}, 
one gets $n>3f^4$, whence the same contradiction emerges.
\eep

\bepr{pr39}  There are no  $D(-1)$--quadruples $(1,b,c,d)$ 
with $10^{13} <b<c<d$,  $b^{1.16}\le c<b^{1.2}$,
and  $b$, $c$ given  by~\eqref{ec2}.
\eepr \bep
Assume the contrary.
By the previous result, $k\ge 4$.  For $k=4$ one gets
\[
f^2(16f^4-12f^2+1)^2+1 <  \left(16f^4(2f^2-1)^2+1\right)^{1.2}
< 169f^{10}+1,
\]
which is false for $f>1$.  For $k=5$ one adapts the reasoning 
introduced in the proof of  Proposition~\ref{pr34} to obtain  
$n> 10^{28}$, in contradiction with~\cite[Proposition~4.3]{bcm}.

As $0.16(2k-1) < 2$ yields $k\le 6$, it remains to consider 
the possibility $k=6$. The argument indicated at the end of
the proof of Proposition~\ref{pr35} can be adapted to the 
present context. One finally obtains $n>3f^4$, which is not 
compatible with the existence  of a  $D(-1)$--quadruple
 subject to all 
constraints from the hypothesis of the present proposition.
 \eep

Summing up what has been done in this section and noticing
that condition~\eqref{ecmild} holds if $f$ has no prime
divisor congruent to $1$ modulo $4$, we get the next result. 

\bet{tef}
There are no  $D(-1)$--quadruples $(1,b,c,d)$  with 
$10^{13} <b<c<d$ and  $b$, $c$ given  by~\eqref{ec2}.
In particular, there exists no $D(-1)$--quadruple for 
which the corresponding
$f$ has no prime divisor congruent to $1$ modulo $4$.
\eet

An alternative proof, more familiar to experts in Diophantine
equations, is based on linear forms in logarithms. Here is the
sketch of such a reasoning.

For the rest of the paragraph we put 
\[
\alpha=s+\sqrt{c}, \quad  \beta = r+\sqrt{b}, \quad 
\mathrm{and}  \quad \gamma=\sqrt{\frac{s\sqrt{b}}{r\sqrt{c}}}. 
\]
As in~\cite{het} and~\cite[Section~4]{bcm}, to a putative 
$D(-1)$--quadruple $(1,b,c,d)$ it is associated a linear form 
in logarithms
\[
 \Lambda := m\log \alpha -l \log \beta +\log \gamma.
\]
We put 
\[
\Delta:= (k+1)m-kl.
\]
Then
\[
(k+1)\Lambda=\log (\gamma ^{k+1} \alpha ^\Delta)-l\log (\beta ^{k+1}
\alpha ^{-k}) 
\]
can be considered as a linear form in the logarithms of two 
algebraic numbers. An elementary study shows that one has
\[
 \vert \log (\beta ^{k+1}\alpha ^{-k})\vert \le \frac{k+1}{f^2}
\quad \mathrm{and}  \quad  \vert \Delta \vert \le \frac{2l}{f^2 \log f}.
\]
Then Theorem~\ref{tef} follows from Laurent's estimates on linear
forms in two logarithms given in~\cite{lau} and our computations
which showed that $f>10^7$ and $b>10^{13}$ for each 
$D(-1)$--quadruple.

\medskip 

Each approach has its own advantages over the other. The 
former is more ``human-friendly'' (and consequently longer),
provides insight and has explanatory power, while the latter 
is computer-intensive and therefore shorter yet less enlightning.
Since the former approach involves ideas who proved to be pivotal
for subsequent developments, we decided to expound it extensively.

The idea at the basis of the proof of Theorem~\ref{tef} can be 
succintly stated ``reduce the master equation to a Pellian equation''.
The same paradigm can be applied for an arbitrary $D(-1)$--quadruple.

Write $f=f_1f_2$, with $f_1$ the product of all the prime 
divisors of $f$ which are congruent to $1$ modulo $4$,
multiplicity included. Then in any solution $(r,s)$ to  
\eqref{ecrs} one has
\[
r=f_2u, \quad s=f_2v,  
\]
for some positive integers  $u$, $v$ satisfying
\beq{ecuv} 
u^2+v^2=2fuv+f_1^2 \Longleftrightarrow (v-fu)^2-(f^2-1)u^2=f_1^2.
\eeq

Below is a specialization of Frattini's theorems from~\cite{frat} 
and~\cite{frat1} giving a representation for the nonnegative 
solutions to the equation relevant for us
\beq{ecf111}
W^2-(f^2-1)U^2=f_1^2.
\eeq

\bepr{prfrat} The nonnegative solutions of Eq.~\eqref{ecf111} are
given by
\[
 w +u \sqrt{f^2-1} =\bigl(w_0+u_0\sqrt{f^2-1} \,\bigr) 
 \bigl(f+\sqrt{f^2-1} \,\bigr)^k, \quad k \ge 0,
\]
or
\[
 w +u \sqrt{f^2-1} =\bigl(w_0-u_0\sqrt{f^2-1} \,\bigr) 
 \bigl(f+\sqrt{f^2-1} \,\bigr)^k, \quad k \ge 1,
\]
where $(w_0,u_0)$ runs through the nonnegative solutions of 
Eq.~\eqref{ecf111} with 
\[
 f_1\le w_0 \le f_1\sqrt{\frac{f+1}{2}}, \quad 
 0\le u_0 \le \frac{f_1}{\sqrt{2(f+1)}}.
\]
\eepr

Let $\gamma=f+\sqrt{f^2-1}$ be the fundamental solution to the Pell
equation $X^2-(f^2-1)Y^2=1$ and $\eps=w_0 +\zeta u_0\sqrt{f^2-1}$,
where $\zeta \in \{-1,1\}$, a fundamental solution as described above. 
According to Lemma~3.5 from~\cite{bcm}, $v>3.999^{1/2}f u$, so that
all positive solutions to Eq.~\eqref{ecuv} have the form
\beq{ec4p}
v-fu +u \sqrt{f^2-1}=\eps \gamma ^k, \quad k\ge (1 -\zeta)/2.
\eeq
Introducing the algebraic conjugates $\vga=f-\sqrt{f^2-1}$,
$\veps=w_0 -\zeta u_0\sqrt{f^2-1}$, one readily obtains
\begin{align}
r &= \frac{f_2\bigl(\eps \, \gamma ^k-\veps \, \vga ^k
\bigr)}{\gamma -\vga},  \label{formr} \\
s &= \frac{f_2\bigl(\eps \, \gamma ^{k+1}-\veps \, \vga ^{k+1}
\bigr)}{\gamma -\vga},  \label{forms}
\end{align}
whence
\begin{align}
b &= \frac{f_2^2\bigl(\eps \, \gamma ^k+\veps \, \vga ^k
\bigr)^2-4}{(\gamma -\vga)^2},  \label{formb} \\
c &= \frac{f_2^2\bigl(\eps \, \gamma ^{k+1}+\veps \, \vga ^{k+1}
\bigr)^2-4}{(\gamma -\vga)^2}.  \label{formc}
\end{align}

Note that when $u_0=0$ one has $\eps=\veps=f_1$, so that 
\eqref{formr} and \eqref{forms} coincide with~\eqref{ec1} 
and~\eqref{ec2}, respectively.

One major source of difficulties with this approach is the fact 
that the components $w_0$, $u_0$ of a fundamental solution are 
known only approximately, being confined to a box defined by the 
inequalities stated in the last line of Proposition~\ref{prfrat}.
Another reason for complexity is the existence of positive solutions 
to Eq.~\eqref{ecf111} for which $\zeta=-1$. We have succeeded to 
overcome all such complications and prove Theorem~\ref{tede} along 
these lines. Our attempts to simplify the proof and avoid intricate 
arguments were successful as soon as we changed once more the 
underlying paradigm. 

Multiplication by a power of the minimal solution for the 
associated Pell equation can be viewed as a vehicle to move 
from a fundamental solution to Eq.~\eqref{ecf111}  to a 
solution of interest. Metaphorically speaking, one can say that 
in the proof for Theorem~\ref{tede} presented in Section~\ref{sec4} 
we travel backwards --- we examine to what extent  information 
about a specific solution  is transferred to associated solutions.

Before making explicit the explanations  alluded to above, we 
present in the next section strenghtened versions for some 
technical results from~\cite{bcm}.

\section{Bounds for linear forms in logarithms} \label{secAleks}

Recall that for a nonzero algebraic number $\gamma$ of degree $D$ 
over $\QQ$, with minimal polynomial $A\prod _{j=1}^D(X-\gamma^{(j)})$ 
over $\ZZ$, the absolute logarithmic height is defined by
\[
\hh (\gamma) = \frac{1}{D}\left( \log A +\sum_{j=1}^D \log^+ \mid (\gamma^{(j)}) \mid\right),  
\]
where $\log^+ x=\log \max (x,1)$.

Next we quote from~\cite{alex} a theorem giving very good lower
bounds for linear linear forms in the logarithms of three algebraic
numbers under hypotheses that are easily checked in the context 
of interest here. 

\bet{Alex} \emph{(Aleksentsev)}
Let $\Lambda_1$ be a linear form in logarithms of $n$ multiplicatively 
independent totally real algebraic numbers $\beta_{1}, \ldots \beta_{n}$, 
with rational coefficients $b_{1}, \ldots, b_{n}$. Let $\hh(\beta_{j})$ 
denote the absolute logarithmic height of $\beta_{j}$ for $1\leq j \leq n$. 
Let $D$ be the degree of the number field 
$\mathcal{K} = \QQ (\beta_{1}, \ldots, \beta_{n})$, and let 
$B_{j} = \max(D \hh(\beta_{j}), |\log \beta_{j}|, 1)$. Finally, let
\beq{tail}
E = \max\left( \max_{1 \leq i, j \leq n} \left\{ \frac{|b_{i}|}{B_{j}} + \frac{|b_{j}|}{B_{i}}\right\}, 3\right).
\eeq
Then
\begin{equation*}\label{kidney}
\log |\Lambda_1| \geq - 5{.}3  n^{-n+1/2} (n+1)^{n+1}(n+8)^{2}(n+5)(31{.}44)^{n} 
D^{2} (\log E)\log(3nD) \prod_{j=1}^n  B_{j} .
\end{equation*} 
\eet
We apply Theorem \ref{Alex} for $D=4$, $n=3$, and 
\beq{eqlam1}
 \Lambda_1 = 2m\log\beta_{1} - 2l\log\beta_{2} + \log\beta_{3},
\eeq
with the choices 
\beq{cheek}
\beta_{1}  = s+ \sqrt{c}, \quad \beta_{2}  = r + \sqrt{b}, \quad 
\beta_{3} = \frac{s\sqrt{b}}{r\sqrt{c} },  \quad
b_1  = 2m, \quad b_2  = -2l, \quad b_3  = 1.
\eeq

The required multiplicative independence readily follows by noting that 
$\beta_{1}$ and $\beta_{2}$ are algebraic units while $\beta_{3}$
is not. Indeed, any possible relation of multiplicative dependence 
has the shape $\beta_1^u = \beta_2^v$ for some positive integers $u$, 
$v$. Note that $\QQ (\beta_1) \cap \QQ (\beta_2) =\QQ$, as otherwise
$b$ and $c$ would have the same square-free part, so that $bc$ would 
be a perfect square, in contradiction with $bc-1=t^2$. One concludes
that it holds $\beta_1^u \in \QQ$, which is not possible because 
$\beta_1$ is not a root of unity.

For compatibility with~\cite{bcm}, we introduce the notation
$\alpha=s+\sqrt{c}$, $\beta=r+\sqrt{b}$. It is clear that it holds
 \[
 \hh (\beta_1) =\frac{1}{2} \log \alpha , \quad \hh (\beta_2) =\frac{1}{2} \log \beta ,  
 \]
so that
\beq{eq:valB}
B_1 = 2\log \alpha , \quad B_2= 2\log \beta .
\eeq 
The minimal polynomial for $\beta_3$ is $r^2cX^2-s^2b$ divided by 
$\gcd (r^2c,s^2b)$, so that 
\[
\hh (\beta_3) = \frac{1}{2} \log \left(\frac{s^2 b}{\gcd (r^2c,s^2b)}\right).
\]
As the lower bound for $\log |\Lambda_1|$ given by Aleksentsev's 
theorem decreases when $B_3$ increases, we can take
\[
B_3 =4\log s\sqrt{b}. 
\]
Combining the obvious relations $\alpha > \beta > \beta_3$, $B_3>B_1>B_2$
with $m\log \alpha < l\log \beta$ (proved in Lemma 3.3 from~\cite{het})
and  its consequence $l>m$,  one obtains
\[
 E=\max \left( \frac{2l}{\log \beta}, 3 \right).
\]
Having in view that by Theorem~\ref{tenoi} one has $b<10^{148/1.16}$, 
for $l \ge 250$ one gets
\[
 E= \frac{2l}{\log \beta}.
\]

Since $\log \Lambda_1 < -4l\log \beta + \log (b)-\log (b-1)$ by Lemma~3.1
from~\cite{het}, one has
\[
 l < 6.005171\cdot 10^{11} \log \alpha \log (s\sqrt{b}) \log \left(\frac{2l}{\log \beta} \right).
\]

Most of the previous work on $D(-1)$-quadruples has focused on the $z$-component
of the solutions to system~\eqref{ecpc}--\eqref{ecpb}. In order to
use the information already available in the literature, we shall derive 
from the inequality above one involving $m$ and subsequently another one 
in terms of $n$.  In a first step towards this goal
we employ the elementary fact that the function $x\mapsto x/\log x$ is
increasing for $x>3$. By Lemma~3.3 from~\cite{het}, we thus get
\[
 m < 6.005171\cdot 10^{11} \log \beta \log (s\sqrt{b})
 \log \left(\frac{2m \log \alpha}{(\log \beta )^2} \right).
\]

A slight simplification is possible thanks to  the following result.
 
\bel{lecomp}
$
\displaystyle \frac{\log (s+\sqrt{c})}{\log (r+\sqrt{b})} < \frac{\log c}{\log b}.
$
\eel \bep 
Consider the real functions $ f_1(x)=\log (\sqrt{x}+\sqrt{x-1})$ 
 and $f_2(x)=\log x$  defined for $x\ge 1$.
As $f_2 '(x) >0$ and $f_1 '(x)/f_2'(x)$ is decreasing for  $x>1$,
by~\cite{avv} we know that
\[
 \frac{\log (\sqrt{x}+\sqrt{x-1})}{\log x} = \frac{f_1(x)-f_1(1)}{f_2(x)-f_2(1)}
\]
is decreasing as well.
\eep

Using this observation together with the obvious inequality
$bs^2 <  bc-1 $, we get
\[
 2m < 6.005171\cdot 10^{11} \log \beta \log (bc-1)
 \log \left(\frac{2m \log c}{\log b \log \beta } \right).
\]
As explained above, for $m\ge 250$ one has $2m>3\log \beta$,
so one can apply the same reasoning to pass from $m$ to $n$
with the help of the inequality $m\log (4c) > n\log (bc-1)$
proved in~\cite[Lemma~2.7]{bcm}. The resulting formula is
\[
 2n < 6.005171\cdot 10^{11} \log \beta \log (4c)
 \log \left(\frac{2n \log (bc-1) \log c}{\log b \log (4c)\log \beta } \right).
\]

Since $2r < \beta < 2\sqrt{b} $, Theorem \ref{Alex} yields the following corollary.
\bec{aleup}
If $n\ge 250$, then
\[
 n < 1.5002\cdot 10^{11} \log (4b) \log (4c)
 \log \left(\frac{4n \log (bc)}{\log b \log (4b) } \right).
\]
\eec 

Upper bounds of this type are complemented by reverse inequalities. 
Our next immediate goal is to  sharpen some lower bounds for $n$ 
in terms of $b$ and $c$ established in~\cite{bcm}. To this end,
we shall use the positive integer $A$ introduced in~\cite{dufu} 
by formula  
\[
A=(2b-1)c-2rst. 
\]
Routine calculations lead to the simpler 
statement $A=f^2+b$. For the proof of our next results we recall 
from Lemma~3.4  of~\cite{bcm} that $A$ satisfies the double inequality
\beq{eqet2}
 \frac{c-5}{4b}+b < A < \frac{1}{3.999} \left(\frac{c}{b} +4b\right)
\eeq
as well as the congruence
\beq{congA}
2(bn^2-m^2) \equiv \pm An \pmod c. 
\eeq 
Occasionally we shall rewrite this as
\beq{eqA}
2(bn^2-m^2) +\rho An = jc
\eeq 
for a fixed $\rho\in \{-1,1\}$ and a certain integer $j$.

Actually, slightly stronger upper bounds on $A$ are valid in the 
context of interest in this paper.

\bel{letarA}
Let $(1,b,c,d)$ with $1<b<c<d$ be a  $D(-1)$--quadruple. 
Then:

\emph{a)} For $c<b^3$ one  has
\[
A<\left(\frac{c}{ 4b}+b\right) \left(1+\frac{1 }{ b}\right).
\]

\emph{b)} $A < 2b$ for $c< 4\, b^2$.

\emph{c)} $A< c/(3.9999\, b)$ for $c>  b^3$.
\eel \bep
For part a) we use the inequality 
\[
A < b + \frac{1 }{ 4}\left( \frac{b-1}{ c-1}+\frac{c-1}{ b-1}\right) +\frac{1 }{ 2}
\]
established in the proof of Lemma~3.4 from~\cite{bcm}.
Since
\[
\frac{b-1}{ c-1} < \frac{b}{ c}
\]
because $b<c$, and
\[
\frac{c-1}{ b-1} < \frac{c}{ b-1}= \frac{c }{ b}+\frac{c }{ b^2}+\frac{c}{ b^2(b-1) },
\]
the result follows from the hypothesis $c<b^3$ and the estimate $b>10^{13}$.

 When  $c< 4\, b^2$, part a) yields $A< 2b+2$. The assumption $c<4b^2$
implies $s <2\, b$, so for $s=2b-1$ one sees from the definition
of $A$ that $2b-1$ divides  $A$, so one necessarily has $A=2b-1$. 
Hence, $c$ is odd, which is not possible with  $s$ odd. 

When $s=2b-2=2r^2$, one readily gets $t=2r^3+r$ and $f=t-sr=r$.
As we already know that there is no  $D(-1)$--quadruple with 
$\gcd (r,s)=f$, we conclude that $s\le 2b-3$. Therefore,
\[
 A< \frac{c}{4b} +b +2 \le 2b -1+\frac{9}{4b} <2b,
\]
as claimed in b).

The inequality  c) follows as soon as we show that it holds
\[
 b+1+\frac{c-1}{4( b-1)} < \frac{c}{3.9999\, b}.
\]
This is a corollary of the slightly stronger inequality
\[
 b+1 < \frac{bc-20000 \, c}{15996\, b(b-1)}
\]
valid because $c>b^3$ and $b>10^{13}$.
\eep

Now we are in a position  to give a simplified list of 
lower bounds for $n$ in terms of $b$ and $c$. More precisely,
the constants appearing in these bounds improve upon those
provided by Lemmata~3.6,~3.9, and~4.2 from~\cite{bcm}.

\bepr{prmarg}
Let $(1,b,c,d)$ with $1<b<c<d$ be a  $D(-1)$--quadruple
with $c >  b^3$. Then $n > \min \{b,0.125 \, \sqrt{c/b} \}$. 
\eepr \bep
We reason by reduction to absurd. Assuming that the thesis is
false, with the help of Lemma~\ref{letarA} we get
\[
 0 < An \le \frac{ c}{3.9999} < 0.2501 \, c,
\]
\[
 0 < 2(bn^2-m^2) <2bn^2  \le \frac{ c}{32} < 0.0313\, c.
\]
Therefore, congruence~\eqref{congA} is actually an equality 
$An=2(bn^2-m^2)$, whence $A < 2bn \le 0.25\, \sqrt{bc}$. This 
inequality is not compatible with $A > c/(4b) $ when $c>b^3$.
\eep

\bepr{prmarg2}  
Let $(1,b,c,d)$ with $1<b<c<d$ be a  $D(-1)$--quadruple
with $4b^2 < c <  b^3$. Then $n > 0.125 \, c/b^2$. 
\eepr \bep
As before, we supppose that the conclusion is false. Note
that part a) of Lemma~\ref{letarA} entails $A<0.6 \, c/b$,
whence $An < 0.075 \, c^2/b^3$. Since $2bn^2 < 0.032 \, c^2/b^3$,
we conclude yet again that congruence~\eqref{congA} is actually
an equality. We  therefore get 
\[
 \frac{c}{4b} < A<2 \, bn \le \frac{c}{4b},
\]
a blatant contradiction.
\eep

The result just proved is not useful when $c$ is close to $4b^2$.
One way to eliminate this inconvenience follows.   In the statement below 
we refer to Eq.~\eqref{eqA}.

\bepr{pr3.8} Suppose $(1,b,c,d)$ with $1<b<c<d$ is a  $D(-1)$--quadruple
with $ c <  b^3$.

\emph{a)} If $\rho =1$, then $n>0.5 \sqrt{c/b}$.

\emph{b)} Let $\rho =-1$. Then $j$ is nonnegative. If $j$  is positive,
then $n>0.5 \sqrt{2c/b}$. If $j=0$, then $c> 7164532\, b^2 > b^{2.155}$ and 
\[
 n> \left\{\begin{array}{ll}
     c^{2/11} & for \ c\ge \max \{ b^{2.5}, 10^{50} \}, \\
     0.214\, (c/b)^{1/3} & for \  c< b^{2.5}.
    \end{array}
    \right. 
\]
\eepr  \bep 
The result is very close to Lemma~3.8 from \cite{bcm}. There are
two differences: the hypothesis $c<b^3$ (instead of $c<b^{2.75}$)
which allows one to employ part a) of the above Lemma~\ref{letarA}
and the conclusion $c> 7164532\, b^2$ (instead of $c> 51.99\, b^2$).

a) The proof given in~ \cite[Lemma~3.8]{bcm} is valid under the 
present hypotheses.

b) In \emph{loc. cit.} it was shown that for $\rho =-1$ and $j=0$ one has
$n>0.214\, (c/b)^{1/3}$ when $ c< b^{2.5}$. Since $b>20^{10}$, one gets 
\[
 A = 2bn -2m^2/n >2(b-4)n > 0.4279 \bigl(51.99\cdot 20^{10} \bigr)^{1/3} b.
\]
As a consequence of  Lemma~\ref{letarA} a) one also has
\[
 A<1.0001 \left(\frac{c}{4b}+b \right).
\]
Comparison of the two bounds on $A$ results in the inequality $c>138703\, b^2$.

Resume the reasoning from the previous paragraph with this information
instead of $c> 51.99\, b^2$. The outcome is an improved lower bound on $c$.
After fourteen more iterations one obtains  $c> 7164532\, b^2$. According
to~\cite[Proposition~4.1]{bcm}, for $b^2 <c<b^3$ one has $b<10^{44}$.
This readily gives $7164532 > b^{0.155}$, which ends the proof.
\eep 

For hypothetical $D(-1)$--quadruples with $c<4\, b^2$ we also offer two 
kinds of lower bounds for $n$.

\bepr{pr3.1}
Let $(1,b,c,d)$ with $1<b<c<d$ be a  $D(-1)$--quadruple. 
If  $ c < 4\,  b^{2}$, then   $n > 0.707 \, \sqrt{c/b}$.
\eepr \bep 
The previous proposition gives $n> 0.5 \sqrt{2c/b}$ when $\rho=-1$, 
which is slightly stronger than the claimed inequality. When 
$\rho=1$ then one has  $2(bn^2-m^2)+An \ge c$. In view of our 
previous Lemma~\ref{letarA}, we get $2bn^2+2b n-c > 0$, so that
\[
 2bn > -b +\sqrt{b^2+2bc}> 1.414 \, \sqrt{bc},
\]
where the last inequality holds because $c>3.999 \, f^2 b > 3.999 
\cdot 10^{14}b$ by Lemma~3.5 of~\cite{bcm}.
\eep 

\bepr{pr3.9} Suppose $n\ge 1000$. 

\emph{i)} If $b^{1.32} \le c <b^{1.40}$, then 
 $\displaystyle  n \ge  \left(  \frac{15.927 \, b^2}{ c}\right)^{1/4}$ .
 
\emph{ii)} If $b^{1.27} \le c <b^{1.32}$, then 
 $\displaystyle  n \ge  \left(  \frac{15.830 \, b^2}{ c}\right)^{1/4}$ .

\emph{iii)} If $b^{1.22} \le c <b^{1.27}$, then 
 $\displaystyle  n \ge  \left(  \frac{15.387 \, b^2}{ c}\right)^{1/4}$ .

\emph{iv)} If $b^{1.16} \le c <b^{1.22}$, then 
 $\displaystyle  n \ge  \left(  \frac{12.850 \, b^2}{ c}\right)^{1/4}$ .
\eepr \bep 
  The new idea is to use the observation that for any 
integers $L$, $M$, $R$, from $L\equiv R\pmod{2M}$ it follows 
$L^2\equiv R^2\pmod{4M}$.

We put $e=f^2-1$, $\Delta=f^2$, then $A=b+\Delta$ and, as seen 
in the proof of Lemma~\ref{letarA}, it follows  that
\[
\Delta < \frac{c}{ 4b} + \frac{c}{4 b^2} + 1 < 0.2501\, \frac{c }{ b}.
\]
Notice that $b +e-1 = 2ft -c $, so that
\[
b^2+2(e-1)b +e^2-2e+1 \equiv c^2+4f^2t^2 \equiv c^2-4f^2=c^2-4(e+1) \pmod{4c}.
\]
It follows that
\[
b^2+2(e-1)b +e^2+2e+5 \equiv c^2  \pmod{4c}.
\]

The congruence method introduced in~\cite{dup} is based on the relation
$s(m^2-bn^2) \equiv \rho rtn \pmod{4c}$. Multiplying both sides by $2s$
one obtains
\[
2(bn^2-m^2) \equiv -\rho A n +cn \pmod {2c},
\]
equivalently
\[
b(2n^2+\rho n)  \equiv 2m^2-\rho \Delta n +cn \pmod {2c}.
\]
From this we get
\[
b^2 (2n^2+\rho n)^2 \equiv (2m^2-\rho \Delta n +cn)^2 \pmod {4c}
\] 
as well as
\[
2(e-1)b(2n^2+\rho n)^2  \equiv 2(e-1)(2n^2+\rho n) ( 2m^2-\rho \Delta n +cn) \pmod {4c}.
\]
By summation we get that $\bigl(b^2+2(e-1)b \bigr)(2n^2+\rho n)^2$
is congruent modulo $4c$ to
\[
( 2m^2-\rho \Delta n +cn)^2 +2(e-1)(2n^2+\rho n)(2m^2 - \rho \Delta  n + cn),
\]
equivalently, again modulo $4c$,
\[
(c^2-e^2-2e-5)  (2n^2+\rho n)^2 \equiv
 ( 2m^2-\rho \Delta n +cn)^2+ 2(e-1)(2n^2+ \rho n)
(2m^2 - \rho \Delta  n + cn).
\]

Considering separately  even and odd values of $n$, it is seen
\[
- (e^2+2e+5)  (2n^2+\rho n)^2 \equiv
 ( 2m^2-\rho \Delta n )^2
+2(e-1)(2n^2+ \rho n)(2m^2 - \rho \Delta  n ) \pmod {4c}.
\]

Now we want to find an upper bound for the expression
\[
\Phi :=  (e^2+2e+5)  (2n^2+\rho n)^2 +
 ( 2m^2-\rho \Delta n )^2
+2(e-1)(2n^2+ \rho n)) |2m^2 - \rho \Delta n |  .
\]
We proceed piece by piece, taking into account the relative
size of $b$ and $c$ as well as the upper bound on $c$. We give 
all details for part i), leaving the other cases to the reader.

According to Lemma 2.5 from~\cite{bcm}, we have
\[
 \frac{2m-1}{2n} < \frac{\log \bigl(4\, c^{1+1/1.32} \bigr) }{\log \bigl(3.996\, c \bigr)}
 < \frac{\log \bigl(4\cdot 10^{147.43\cdot 58/33} \bigr)}{\log \bigl(3.996\cdot 10^{147.43}\bigr)}
 < 1.7545,
\]
whence 
\[
m< 1.7545 \, n+ 1/2 \le 1.7555 \, n.
\]

Therefore,
\begin{align*}
 |2m^2 - \rho \Delta n | & < 2 \times 1.7555^2 n^2 + 0.2501 \frac{cn }{ b} 
 = \left( 6.16005   +  0.2501 \frac{c}{ nb} \right) n^2, \\
 |2m^2 - \rho \Delta n |^2 & <  \left( 6.16005  +  0.2501 \frac{c}{ nb} \right)^2  n^4,
\end{align*}
\[
2(e-1)(2n^2+\rho n) |2m^2-\rho \Delta n| <  \frac{0.5002 c }{ b}\times 
2.001 \,n^2 \times  \left( 6.16005  +  0.2501 \frac{c}{ nb} \right) n^2,
\]
hence
\[
2(e-1)(2n^2+\rho n) |2m^2-\rho \Delta n|  < 1.001  \frac{c}{ b} 
\left( 6.16005  +  0.2501 \frac{c}{ nb} \right) n^4.
\]
By
\begin{align*}
 e^2+2e+5  =\Delta^2+4 & < \left( 0.2501^2 + \frac{4}{b^{0.64}} \right)\frac{c^2}{b^2}\\
& < \left( 0.2501^2 + \frac{4}{10^{13\cdot 0.64}} \right)\frac{c^2}{b^2}
 < 0.062551 \, \frac{c^2}{b^2}
\end{align*}
we also have
\[
(e^2+2e+5 )  (2n^2+\rho n)^2 < 0.062551 \times  \frac{c^2 }{ b^2} \times 2.001^2 n^4 
< 0.250455 \, n^4 \,  \frac{c^2 }{ b^2} .
\]
Collecting all these  estimates we get that $\Phi $ is less than
\[
\left( 0.250455   +  \left( 6.16005  \frac{b}{ c} +   \frac{0.2501}{ n} \right)^2  
 + 1.001  \left( 6.16005  \frac{b}{ c} +   \frac{0.2501}{ n} \right) \right) \, n^4 \,  \frac{c^2 }{ b^2}.
\]
 
\medskip

 Using the  known lower bounds on $n$ and $c/b$, it is easy to verify that
\[
\Phi   < 0.251134\, n^4 \,  \frac{c^2 }{ b^2} 
\]
and we see that
\[
n <  \left( \frac{15.927 \, b^2}{ c}\right)^{1/4} \ \Longrightarrow\  \Phi < 4\, c.
\]
But, by a previous congruence, the nonnegative integer $\Phi$ is a
multiple of $4c$, so for $\Phi <  4\, c $ it holds
\[   
  (e^2+2e+5)  (2n^2+\rho n)^2 +
 ( 2m^2-\rho \Delta n )^2+2(e-1)(2n^2+ \rho n))(2m^2 - \rho \Delta  n )=0.
\]
The left hand side of the last equation is of the form
\[
(e^2+2e+5) X^2 +2(e-1) XY+Y^2,  \quad {\rm with}\ X = (2n^2+\rho n), 
\ Y =  ( 2m^2-\rho \Delta n ),
\]
a quadratic form whose discriminant (equal to $- 4\, f^2$) is negative, so 
that $\Phi$ is always positive when $n>0$, a contradiction which implies
$n \ge  \left( {15.927 \, b^2}{ c}\right)^{1/4}$.
\eep 

The results just proved serve to improve Theorem~\ref{tenoi}. 
To this end we combine them with Aleksentsev's theorem  
 in conjunction with a similar result due to Matveev~\cite[Theorem 2.1]{Mat98}
applicable in the following context. 

\smallskip 

Let $\beta_1,\beta_2,\beta_3$ be real algebraic 
numbers and denote $K:=\QQ(\beta_1,\beta_2,\beta_3)$. Put $D:=[K:\QQ]$. 
Assume that $\beta_1,\beta_2,\beta_3$ satisfy the Kummer condition, that is, 
\[
[K(\sqrt{\beta_1},\sqrt{\beta_2},\sqrt{\beta_3}):K]=8.
\]
Consider a linear form 
$
\Lambda_1:=b_1 \log \beta_1+b_2 \log \beta_2+ b_3 \log \beta_3,
$
where $b_1,\,b_2,\,b_3$ are integers with $b_3 \ne 0$. 
Put $A_j:=h(\beta_j)$ for $1 \le j \le 3$. 
We take $E,\,E_1,\,C_3,\,C_1,\,C_2$ as follows:
\begin{align*}
E & \ge \frac{1}{3D} \max \left\{\left|\pm \frac{\log \beta_1}{A_1}\pm 
\frac{\log \beta_2}{A_2} \pm \frac{\log \beta_3}{A_3}\right|\right\},\\
E_1&=\frac{1}{2D}\left(\frac{1}{A_1}+\frac{1}{A_2}+\frac{1}{A_3}\right),\\
C_3^*\exp(C_3^*)\frac{Ee}{2}&\ge e^3,\quad C_3=\max\{C_3^*,3\},\\
C_1&=\left(1+\frac{e^{-6}}{148}\right)(3 \log 2 + 2) \frac{4}{3C_3},\\
C_2&=16\left(6+\frac{5}{3 \log 2+2}\right)\frac{e^6}{3^{1/2}C_3}.
\end{align*}
We also put
\begin{align*}
\Omega&:=A_1A_2A_3,\\
\omega&:=\Omega \left(\frac{DC_1}{e}\right)^3C_3 \exp(C_3)\frac{Ee}{2}.
\end{align*}
Let $C_0$ be a real number satisfying
\[
C_0 \ge \max\left\{2C_3,\log\left(4C_2 \max
\left\{\frac{C_0\omega}{4C_1A_3},C_0,\frac{2E_1C_3}{C_1}\right\}\right)\right\}.
\]
Furthermore, put
\begin{align*}
B_0&:=\sum_{j=1}^2 \frac{|b_3|+|b_j|}{8\gcd(b_j,b_3)C_0C_2\omega},\\
B_1&:=\sum_{j=1}^2 \frac{1}{24\gcd(b_j,b_3)C_1}\left(\frac{|b_3|}{A_j}+\frac{|b_j|}{A_3}\right),\\
B_2&=\sum_{j=1}^2 \frac{|\log \beta_j|(|b_3|+|b_j|)}{8|b_3|C_0C_2\omega},\\
B_3&=\sum_{j=1}^2 \frac{|\log \beta_j|}{24|b_3|C_1}\left(\frac{|b_3|}{A_j}+\frac{|b_j|}{A_3}\right),
\end{align*}
and take a real number $W_0$ satisfying
\[
W_0 \ge \max\{2C_3,\log(e(1+B_0+B_1+B_2+B_3))\}.
\]
Now we are ready to state \cite[Theorem 2.1]{Mat98} in a form applicable to our situation. 
\bet{thm:Mat98} \emph{(Matveev)}
Suppose that 
\begin{align*}
2\omega\min\{C_0,W_0\} &\ge C_3,\\
\omega \min\{C_0,W_0\} & \ge 2C_1C_3 \max \{A_1, A_2,A_3\},\\
3(4C_1)^2 4C_0\Omega  & \ge C_3\max \{A_1, A_2,A_3\}.
\end{align*}
Then, 
\[
\log |\Lambda_1| > -11648C_2C_0W_0 \omega.
\]
\eet

In order to apply this result to  the linear form~\eqref{eqlam1},
we need to check the Kummer condition is valid.

First we show that $\sqrt{\beta_1} \not \in K$. 
Assume on the contrary that $\sqrt{\beta_1} \in K$. 
Then one may write
$
\sqrt{\beta_1} =l_0+l_1\sqrt{b}+l_2\sqrt{c}+l_3\sqrt{bc}
$
with $l_0,\,l_1,\,l_2,\,l_3 \in \QQ$. 
Squaring both sides yields
\begin{align*}
s+\sqrt{c}&=l_0^2+bl_1^2+cl_2^2+bcl_3^2+2(l_0l_1+cl_2l_3)\sqrt{b}\\
           &\quad +2(l_0l_2+bl_1l_3)\sqrt{c}+2(l_0l_3+al_1l_2)\sqrt{bc},
\end{align*}
whence
$
s=l_0^2+bl_1^2+cl_2^2+bcl_3^2, \quad  1=2(l_0l_2+bl_1l_3). 
$
The arithmetic mean -- geometric mean inequality yields
\[
 s\ge 2|l_0l_2|\sqrt{c} +2|l_1l_3| b\sqrt{c}\ge 2(l_0l_2+bl_1l_3)\sqrt{c} 
 = \sqrt{c} > s,
\]
a contradiction.

Similarly one proves that $\sqrt{\beta_2} \not \in K$.
To check $\sqrt{\beta_3} \not \in K$, we suppose the contrary and 
get $0=l_0^2+bl_1^2+cl_2^2+bcl_3^2$ for some $l_0,\,l_1,\,l_2,\,l_3 \in \QQ$.
Since $b$, $c>0$, it follows that all $l_j$ are zero, so $\beta_3 =0$,
absurd.

Secondly, assume that $\sqrt{\beta_1} \in K(\sqrt{\beta_2})$. 
Then one may write $\sqrt{\beta_1}=k_0+k_1\sqrt{\beta_2}$ for 
some $k_0,\,k_1 \in K$, equivalently 
$
s+\sqrt{c}=k_0^2+k_1^2(r+\sqrt{b})+2k_0k_1\sqrt{\beta_2}.
$
If $k_0k_1 \ne 0$, then this equation shows that $\sqrt{\beta_2} \in K$, 
which is impossible as seen above. If $k_1=0$, then $s+\sqrt{c}=k_0^2$, 
which contradicts $\sqrt{\beta_1} \not \in K$.  It remains $k_0=0$, so that 
$s+\sqrt{c}=k_1^2(r+\sqrt{b})$ with 
$k_1=l_0+l_1\sqrt{b}+l_2\sqrt{c}+l_3\sqrt{bc}$ and 
$l_0,\,l_1,\,l_2,\,l_3 \in \QQ$.  Identification of coefficients of 
$\sqrt{b}$ on the two sides of this equation followed by aplication
of the arithmetic mean -- geometric mean inequality results in
\begin{align*}
0 & = l_0^2+bl_1^2+cl_2^2+bcl_3^2 +2(l_0l_1+cl_2l_3)r \\
 & \ge 2 \bigl(|l_0l_1| + c |l_2l_3|\bigr)\sqrt{b} +2(l_0l_1+cl_2l_3)r \\
 & \ge 2  \bigl(|l_0l_1| + l_0l_1 +c (|l_2l_3|+ l_2l_3 )\bigr) r \ge 0.
\end{align*}
The middle inequality is strict unless $l_0l_1 = l_2l_3 =0$, in which 
case all $l_j$ are zero. In either case we reached a contradiction.
Similarly one shows that $\sqrt{\beta_2} \not \in K(\sqrt{\beta_1})$. 

It remains only to show that $\sqrt{\beta_3} \not \in K(\sqrt{\beta_1},\sqrt{\beta_2})$. 
Assume  the contrary and put
\[
\sqrt{\beta_3}=k_0+k_1\sqrt{\beta_1}+k_2\sqrt{\beta_2}+k_3\sqrt{\beta_1 \beta_2}
\]
with some $k_0,\,k_1,\,k_2,\,k_3 \in K$. 
Squaring both sides, one has
\begin{align}\label{chi'2}
\beta_3 &=k_0^2+k_1^2 \beta_1+k_2^2 \beta_2+k_3^2 \beta_1 \beta_2+2(k_0k_1+k_2k_3 \beta_2)\sqrt{\beta_1}\\
     &\quad +2(k_0k_2+k_1k_3 \beta_1)\sqrt{\beta_2}+2(k_0k_3+k_1k_2)\sqrt{\beta_1 \beta_2}. \notag
\end{align}
If  $k_0k_1+k_2k_3\beta_2\ne 0$, with the help of $\sqrt{\beta_2} \not \in K$
one deduces first that $k_0k_1+k_2k_3\beta_2 + 2(k_0k_3+k_1k_2)\sqrt{\beta_2}\ne 0$
and next that  $\sqrt{\beta_1} \in K(\sqrt{\beta_2})$. A similar 
contradiction is reached assuming either $k_0k_3+k_1k_2 \ne 0$ or 
$k_0k_2+k_1k_3 \beta_1 \ne 0$. So it holds
\[
 k_0k_1= -k_2k_3 \beta_2, \quad   k_0k_2= -k_1k_3\beta_1, \quad k_0k_3= -k_1k_2,
\]
whence
\[
 k_0k_1(k_2^2-\beta_1 k_3^2)=0, \quad k_0k_2(k_1^2-\beta_2 k_3^2)=0, \quad
 k_1k_2(k_0^2-\beta_1 \beta_2 k_3^2)=0.
\]
Having in view what we already proved, it is readily seen that the last 
three equations imply that precisely one of $k_j$ is nonzero. Note that
$k_0 \ne 0$ gives $\sqrt{\beta_3} \in K$, absurd. For $k_1 \ne 0$
one has $s\sqrt{b} = k_1^2 (s+\sqrt{c}) r \sqrt{c}$. Passing to $\QQ$
and comparing the coefficients of $1$ and $\sqrt{c}$ in both hand sides, 
one gets a linear system of equations $sX+cY=0$, $X+sY=0$, with 
$X=l_0^2+bl_1^2+cl_2^2+bcl_3^2$, $Y=2(l_0l_2+bl_1l_3)$, and 
$l_0,\,l_1,\,l_2,\,l_3 \in \QQ$. Since the determinant of this system
is $s^2-c=-1$, it has only the trivial solution, which gives the 
contradiction $k_1=0$. Similarly one can conclude that neither 
$k_2 \ne 0$ nor $k_3 \ne 0$ is possible. 

\smallskip

Now the verification that Kummer condition holds for our $\Lambda$
is complete, so we can proceed with choosing suitable values for the 
parameters in the statement of Theorem~\ref{thm:Mat98}.

As discussed in conection with Theorem~\ref{Alex}, we  take
\[
 A_1 = \frac{1}{2}\log (s+\sqrt{c}) , \quad A_2= \frac{1}{2}\log (r+\sqrt{b}),
\quad A_3 =\log (s\sqrt{b}).
\]
Then
\[
 E\ge \frac{1}{12}\max \left\{ \left \vert\pm 2 \pm 2 \pm 
 \frac{\log (s\sqrt{b}/r\sqrt{c})}{\log (s\sqrt{b})}\right\vert \right\}.
\]
From 
\[
 \frac{s\sqrt{b}}{r\sqrt{c}} =\sqrt{1+\frac{s^2-r^2}{(s^2+1)r^2}}
 <\sqrt{1+\frac{1}{r^2}} < 1+\frac{1}{2r^2},
\]
we see that we can take 
\[
 E=\frac{4+10^{-15}}{12}.
\]
Thus, we may take $C_3^*=2.8$ and $C_3=3$. 

In order to fix a value for $E_1$, we need lower bounds 
for $\log \beta _j$. Using Theorem~\ref{tenoi}, it is readily
seen that a suitable value is 
\[
 E_1=0.033653.
\]

It is easy to see that $C_0$ should satisfy
\begin{align*}
C_0 \ge \log \left(\frac{C_0C_2\omega}{C_1A_3}\right)
    =\log(C_0T)
\end{align*}
with $T=96Ee\,C_1^2 C_2 A_1A_2$, which allows us to take 
\[
C_0=\log T+\log(\log T)+\log(\log(\log T))+2\log(\log(\log(\log T)))
\]
(note that $\log(\log(\log(\log T)))>0$). 
Since $m \log \beta_1 < l \log \beta_2$ and $A_1 > A_2$, one has
\begin{align*}
B_0+B_1+B_2+B_3&<\left(\frac{l}{2C_0C_2\omega}+\frac{1}{6C_1}
\left(\frac{1}{\log \beta_2}+\frac{l}{A_3}\right)\right)(1+\log \beta_1).
\end{align*}
We therefore take
\[
W_0=1+ \log\left( 1+\left( \frac{l}{2C_0C_2\omega}+\frac{1}{6C_1}
\left(\frac{1}{\log \beta_2}+\frac{l}{A_3}\right)\right)(1+\log \beta_1)
\right).
\]
Hence, combining the estimate in Theorem~\ref{thm:Mat98} with 
$0 < \Lambda_1 < \bigl(8ac/(b-1)\bigr)\beta_2^{-4l}$ one gets 
\beq{mal}
l < 69888 C_0 C_1^3 C_2 W_0 E e \log (s+\sqrt{c}) \log (s\sqrt{b})
+0.25 \log \left(\frac{b}{b-1} \right).
\eeq

As previously did, we pass from this inequality to one involving
$m$ and subsequently to one in terms of $n$. Assuming that it holds
$b^{del} < c < b^{Del} $ for some real numbers $1.16< del < Del <4.1$,
we finally get
\beq{ineqM98n}
n < 17472 C_0 C_1^3 C_2 W_0 S e \log (4b) \log (4c),
\eeq
with
\[
 S=1+ \log\left( 1+\left( (Del+1) \left(\frac{1}{2C_0C_2\omega}+
 \frac{1}{6C_1A_3} \right) n + \frac{1}{6C_1\log \beta_2}\right)
 (1+\log \beta_2) \right).
\]

At this moment we have all ingredients for the proof of the main 
result of this section.

\bet{prale}
Let $(1,b,c,d)$ with $1<b<c<d$ be a  $D(-1)$--quadruple. 
Then $b>1.024\cdot 10^{13}$ and $ \max \{10^{14}b,  b^{1.233} \}
 < c< \min \{  b^{2.93}, 10^{99} \}$. 
More precisely:
\begin{enumerate}
\item [i)] If $ b^{2}\le c <  b^{2.93}$, then 
            $b< 6.89\cdot 10^{32}$ and $c< 8.48\cdot 10^{70}$.
\item [ii)] If  $b^{1.5}\le c<  b^2$, then 
            $b< 1.26\cdot 10^{49}$ and  $c <  4.48\cdot 10^{73}$.
\item [iii)] If $b^{1.4} \le c< b^{1.5}$, then 
            $b< 2.07 \cdot 10^{62}$ and  $c <  1.77\cdot 10^{87}$.
\item [iv)] If $b^{1.3} \le c< b^{1.4}$,  then 
            $b< 6.26 \cdot 10^{73}$ and  $c <   10^{99}$. 
\item [v)] If $b^{1.233} \le c< b^{1.3}$, then 
            $b<  10^{69}$ and  $c <  4.85\cdot 10^{89}$. 
\end{enumerate}
\eet 
\bep 
Each interval $b^\gamma < c < b^\delta$ has been covered by subintervals
$b^\mu < c < b^{\mu +0.0001}$. On each subinterval, Corollary~\ref{aleup}
and Proposition~\ref{prmarg} produce an upper bound on $b$, which 
in turns leads to  a bound on $S$. When $n<1000$, instead of 
Proposition~\ref{prmarg} we apply similar results from~\cite{bcm}
valid for $n\ge 7$, which results in much sharper bounds on $b$.
Using the estimate on $S$ in~\eqref{ineqM98n}, an improved upper 
bound on $b$ is obtained. From our computations we learned that
$b<10^{13}$ for $c> b^{2.928}$, whence the conclusion that no
$D(-1)$--quadruple has $c> b^{2.928}$.

We also bound from above $f$ with the help 
of the master equation, which shows that 
\[
f<\frac{s}{2r}+\frac{r}{×2s}.
\]
Since our computations yield that for $c \le b^{1.233}$ one has 
$f< 10^7$, by Proposition~\ref{prexp} we conclude that there 
exists no $D(-1)$-quadruple with $c$ so close to $b$.
\eep 

Comparison with Theorems~\ref{tenoi} reveals superiority of
Theorem~\ref{prale}. However, it is also apparent that a lot
of work is required to confirm the nonexistence of 
$D(-1)$--quadruples  by using tools already employed.
Therefore, completely different ideas are required for further 
advancements. The next section details and clarifies the change
in viewpoint on the problem.

\section{A proof for the main theorem}\label{sec4}

Recall that we have denoted by $J$ the set of pairs of integers 
$(r,s)$ such that there exists a $D(-1)$--quadruple $(1,b,c,d)$ 
with $1<b<c<d$ and $b=r^2+1$, $c=s^2+1$. For $(r,s)\in J$ we put  
$s=r^\theta$ and define
\[
\theta^-=\inf_{(r,s)\in J} \theta, \quad 
\theta^+=\sup_{(r,s)\in J} \theta.
\]
Using the upper bound $c\le 2.5 \, b^6$ (see Theorem~\ref{tenoi}), 
a computer-aided search described in Section~2 of~\cite{bcm} led
to the conclusion that any hypothetical  $D(-1)$--quadruple
satisfies  $r>32\times 10^5$, so that $b>10^{13}$. Hence
\[
  \theta^+ \le 6{.}04.
\]

Based on a refinement of a Diophantine approximation result of 
Rickert~\cite{ric}, Filipin and Fujita proved in~\cite{ff} the 
inequality $c \le 9.6\, b^4$, which yields
\[
  \theta^+ \le 4{.}08.
\]
What we just proved in Theorem~\ref{prale} entails
\[
 \theta^+< 3.
\]

The lower bound $\theta^- > 1{.}16$ has been obtained in~\cite{bcm}
and improved in Theorem~\ref{prale} above to
\[
 \theta^- > 1{.}23. 
\]

Further shortening of the interval $[\theta^-, \theta^+]$  
along these lines becomes hopeless, so we introduce the
new approach mentioned in Introduction. 
Our next concern is to have a closer look at solutions of the 
master equation compatible with the information gathered so far.
A convenient tool was suggested by the fact that any solution 
$(x,y)$ to a Diophantine equation of the type $X^2-2fXY+Y^2=C$ 
gives rise to other two solutions, namely $(-x+2fy,y)$ and 
$(x,2fx-y)$, see~\cite{own}. In this section  we shall see how 
changes of variables indeed allow to transfer  information 
regarding one specific solution  to an associated solution.

Introduce a new variable 
\[F:=s-2rf.
\]
As we shall show shortly, it satisfies
\beq{ecesd}
F \sim  \left\{ \begin{array}{rl}
\frac{1}{4} r^{\theta - 2}, & \mbox{  if} \quad \theta>2,
\\
- r^{2-\theta}, & \mbox{  if} \quad \theta<2.
\end{array} \right.
\eeq
Therefore, we hope to exploit this variable in order to split 
the interval  $[\theta^-,\theta^+]$  into two subintervals 
having a common end-point  about $2$. 

We study $F$ with the help of the equation
\[
r^2+s^2=2frs+f^2
\]
or its equivalent forms
\beq{ec2d}
sF=f^2-r^2,
\eeq
\beq{ec3d}
F^2+2frF+r^2-f^2=0.
\eeq
From the master equation one obtains
\[
f=-rs+(r^2s^2+s^2+r^2)^{1/2}=\frac{s^2+r^2}{rs+(r^2s^2+s^2+r^2)^{1/2}}
\sim \frac{1}{2} r^{\theta-1}.
\]
It follows
\[
F = \frac{f^2-r^2}{s}
\sim \frac{1}{4}r^{\theta-2} -r^{2-\theta},
\]
whence estimate~\eqref{ecesd}. 

The first properties of $F$ are almost obvious.

\bel{le1d}
\emph{a)}  $ \quad F=0 \iff s=2rf \iff f=r  \iff s=2r^2 \iff s=2f^2 $.

\emph{b)}   $ \quad F>0  \iff s> 2rf \iff f>r \iff s>2r^2 \iff s<2f^2 $.

\emph{c)}  $ \quad F<0  \iff s<2rf \iff f<r \iff s<2r^2 \iff s>2f^2 $.
\eel \bep
\emph{b)} To prove ``$ \, F>0 \iff s<2f^2\, $'', notice that from 
Eq.~\eqref{ecrs} one gets 
\[
 r=sf-\sqrt{s^2f^2-s^2+f^2}=\frac{s^2-f^2}{sf+
\sqrt{s^2f^2-s^2+f^2}},
\]
and the last expression is smaller than $s/(2f)$ precisely when
$(s^2-f^2)(s^2-4f^4)<0$.

For ``$\, F>0 \iff s>2r^2\, $'', use $f=-rs+\sqrt{r^2s^2+r^2+s^2}$.
\eep

Observe that there are no $D(-1)$--quadruples for which the
corresponding $F$ is zero.

\bel{le3d}
$F\ne 0$.
\eel \bep
Suppose, by way of contradiction, that the thesis is false.
Since $F=0$ if and only if $s=2f^2$ and $r=f$, we are in a 
situation  we have dealt with  in Section~\ref{sec2}. There 
it was found that this is possible for no $D(-1)$--quadruple.
\eep

From these results it readily follows 
that if $f\ne r$, then $f$ is comparatively far away $r$. The
quantitative expression is given by the next lemma.

\bel{le2d}
\emph{a)} If $ f>r$, then $ f>2r F\ge 2r$.

\emph{b)} If $  f<r $, then  $0>F>-2f r$. 
\eel \bep 
Part a) follows from Eq.~\eqref{ec3d} rewritten as 
$F^2+r^2=f^2-2rfF$. 

b) In view of Lemma~\ref{le1d}, $F$ is negative when  $  f<r $.
The lower bound for $F$ follows from~\eqref{ec3d} rewritten as 
$F^2+2rfF=f^2-r^2$.
\eep

Now we have all ingredients to show that the newly introduced 
variable $F$ indeed serves to separate values of $c$ smaller
than $4b^2$ from those bigger than this threshold.

\bel{le4d}
 $ f<r \iff c<4b^2$. 
\eel \bep
We know that $ f<r$ holds if and only if $s\le 2r^2-1=2b-3$, which 
in turn is equivalent to $c\le 4b^2-12b+10$. Hence, $c<4(b-1)^2$ 
for $ f<r$. To prove the converse implication, note that $c<4b^2$ 
is tantamount to $s\le 2r^2+1$. For $s=2r^2+1$,  Eq.~\eqref{ecrs} 
becomes a quadratic in $f$ without integer roots, having
discriminant  $4r^6+8r^4+6r^2+1=(2r^3+2r)^2+2r^2+1=
(2r^3+2r+1)^2-4r^3+2r^2-4r$. Thus  one has  $s\le 2r^2$,
with equality prohibited by Lemmata~\ref{le1d} and~\ref{le3d}. 
It remains $s<2r^2$, which, according to the last part of
Lemma~\ref{le1d}, means $f<r$.
\eep

Our next result shows that the existence of $D(-1)$--quadruples is not 
compatible with small values of $F$. 

\bepr{prDmic}
There is no $D(-1)$--quadruple  with $|s-2rf| \le 2\cdot 10^6$. 
\eepr \bep 
This claim can be obtained by the following algorithm. 

Start by rewriting Eq.~\eqref{ec3d} in the form
\beq{eqDD}
 X^2  - (F^2+1) Y^2 = F^2,  \quad \mbox{where} \quad X=f-rF , \ Y =r.  
\eeq
For any $F$, one obvious solution is  $(F^2-F+1,F-1)$. A 
conjecture of Dujella predicts that an equation $X^2-(a^2+1)Y^2=a^2$ 
has at most one positive solution with  $0< Y < |a|-1$ (this readily 
implies the nonexistence of $D(-1)$--quadruples, see~\cite{mrw}). 
In~\cite{mrw}, this claim is checked for $|a| < 2^{50}$, so,  for each 
$F$ with absolute value up to $2\cdot 10^6$ we can find at most one 
exceptional solution $(x_0,y_0)$ with $0< y_0  < |F|-1$.

Next we consider solutions $(x,y)$ to Eq.~\eqref{eqDD} associated 
to either the obvious solution or to the exceptional one. Invert 
the relations $x=f-rF$, $y=r$, $F=s-2fr$ to obtain $r=y$, $f=x+rF$,
$s=F+2rf$. Check if the resulting values for $r$, $s$, $f$ satisfy
the necessary conditions $r>20^5$, $r^{1.23} < s < r^3$, $f>10^7$. 

Finally, apply Baker--Davenport lemma for each solution surviving 
the sieving step and produce a contradiction with a known fact.

We use this procedure for $| F| \le  2\cdot 10^6$. For the last step, 
we performed computations with real numbers of 173 decimal digits.
In all cases, the outcome of the reduction step is $n=1$. This
contradicts~\cite[Proposition 2.2]{bcm}, where it was shown that 
for no $D(-1)$--quadruple is $n<7$ possible.
\eep

Now we are in a position to halve the region where $\theta$ is confined.

\bepr{prb2cb3}
There is no $D(-1)$--quadruple with $4\, b^{2} \le c < b^3$.
\eepr \bep
The key new ingredient is the  observation that for $ c < b^3$
one has 
\[
 2bn >A.
\]
Lemma~\ref{le2d} together with Lemma~\ref{le4d} imply 
$A=f^2+b> 4r^2 F^2$. By Proposition~\ref{prDmic}, the right-hand 
side is greater than $15\cdot 10^{12} b$. We have thus obtained
$n> 7\cdot  10^{12}$. However, explicit computations find that
for  $c\ge 4 \, b^{2}$ one has $n<10^{12}$. 

Coming to the proof of the claim, we note that it was 
explicitly established in case $\rho =1$ during the proof of
Lemma~3.8 from~\cite{bcm}. It is also shown there that when 
$\rho =-1$ then $j\ge 0$, so Eq.~\eqref{eqA} gives 
$2(bn^2-m^2) \ge An$, whence $2bn >A$.
\eep

Further compression of the interval $[\theta^-,\theta^+]$ is 
possible by examining other solutions of the master equation.
 Let us define a recurrent sequence by the relation
 \[
 F_{i+1} = 2 fF_i-F_{i-1}, \quad i\ge 0, \quad F_{-1}=-s,
 \quad F_0=-r.
 \]
 It is readily seen that $F_1=F$ and for any $i\ge -1$ it holds
\beq{eqF} 
F_i^2-2fF_iF_{i+1}+F_{i+1}^2=f^2. 
\eeq
Moreover, for $i\ge 0$ one has 
\beq{eqFP}
 F_i= P_{i-1}s -P_i r,
\eeq  
where $P_{-1}=0$, $P_0=1$, and $P_{i+1}=2fP_{i}-P_{i-1}$ for any 
nonnegative $i$.

As we shall see shortly, all terms of the sequence $(F_i)_{i\ge 1}$
have properties similar to those established above for $F=F_1$.
First we argue that all $F_i$ are nonzero. In view of Theorem~\ref{tef},
it is sufficient to prove the next result.

\bel{lenz}
Assume $F_i=0$ for some $i\ge 1$. Then $r=P_{i-1}f$ and $s=P_{i}f$.
\eel \bep 
Since $F_i=0$ is tantamount to $sP_{i-1}=rP_{i}$, Eq.~\eqref{ecrs}
becomes 
\beq{eq9rev}
(P_{i}^2+P_{i-1}^2-2fP_{i-1}P_{i})r^2=P_{i-1}^2f^2.
\eeq 
As the  expression within parantheses is  
\[
P_{i}^2+P_{i-1}^2-2fP_{i-1}P_{i}= P_{i}^2-P_{i-1}P_{i+1}=P_{0}^2-P_{-1}P_{1}=1, 
\]
the equality~\eqref{eq9rev} implies $r=P_{i-1}f$, whence $s=P_{i}f$. 
\eep 

Note that from $F_i^2-F_{i-1} F_{i+1}=f^2$, for $i\ge 1$ one gets by induction
\[
 F_{i} \sim \frac{1}{4}r^{i\theta-i-1} -r^{i+1-i\theta},
\]
so that 
\[
 F_{i} \sim \left\{ \begin{array}{rl}
\frac{1}{4} r^{i\theta - i-1}, & \mbox{  if} \quad \theta>(i+1)/i,
\\
- r^{i+1-i\theta}, & \mbox{  if} \quad \theta<(i+1)/i.
\end{array} \right.
\]

A reasoning similar to the proof of Lemma~\ref{le2d} yields
the following result.

\bel{le2f}
Assume $i\ge 0$. If $F_iF_{i-1} <0$, then $f>-2F_iF_{i-1}$.
If $F_{i-1}F_{i+1} >0$, then  $F_i^2=F_{i-1}F_{i+1} +f^{2}>F_{i-1}F_{i+1} +10^{14}$.
\eel \bep 
The desired inequalities are obtained by rewriting the master
equation in the equivalent forms $F_i^2+F_{i-1}^2 =f^2+2fF_iF_{i-1}$
and $F_i^2-F_{i-1}F_{i+1}=f^2$ and taking into account Proposition~\ref{prexp}.
\eep

\bepr{prFmic}
For any $D(-1)$--quadruple it holds $|F_i| >  2\cdot 10^6$  for $2\le i\le 5$. 
\eepr \bep 
The algorithm described in the proof of Proposition~\ref{prDmic}
can be adapted for the present context. One necessary modification
is in the second step, now the inversion of the equations 
$x=f+F_{i-1}F_i$, $y=F_{i-1}$, $F_i=P_{i-1}s -P_i r$ gives
$r=P_{i-2} F_i -P_{i-1}y$, $s=P_{i-1} F_i -P_{i}y$, $f=x-F_iy$,
because $P_{i-1}^2-P_{i-2}P_i=1$. The resulting values have to
satisfy $r^{1.23} < s < r^{2.05}$ by Theorem~\ref{prale} and
Proposition~\ref{prb2cb3}.
\eep 

The last two results have the following consequence.

\bec{coFpo}
If $c<4b^2$, then $F_i< -10^7$ for $-1 \le i\le 5$.
\eec \bep 
For $F_{-1}=-s< -r^{1.23}$ and $F_0=-r$, the desired conclusion follows
from Proposition~\ref{prexp} in conjunction with Lemma~\ref{le1d} c).
Explicit computations show that for any hypothetical $D(-1)$--quadruple 
one has $f<3\cdot 10^{12}$. Together with Lemma~\ref{le2f} and 
Proposition~\ref{prFmic}, this implies $F_i < -10^7$ for $1\le i \le 5$. 
\eep

A last new ingredient in the proof of our  main result  is 
obtained by applying a specialization of the binomial theorem
\[
0< u< 1 \Longrightarrow \sqrt{1+u} =1+\frac{1}{2}u -\frac{1}{8} u^2 + \frac{L}{16}u^3,
 \quad \mathrm{where}   \quad  0<L<1,
\]
to formula
\[
 f=rs\left(-1+\sqrt{1+r^{-2} +s^{-2}} \right).
\]

Maybe it is worth mentioning that $L$ depends on $u$ but the only
property of the function $L$ used below is its  boundedness.

When one uses the resulting expression
\[
 f=\frac{s}{r} +\frac{r}{s} -\frac{s}{4r^3}-\frac{1}{2rs} -\frac{r}{4s^3}
 +\left( \frac{s}{8r^5}+\frac{3}{8r^3s} +\frac{3}{8rs^3} +  \frac{r}{8s^5}\right) L
\]
in $F_1=s-2fr$, it gives
\[
 F_1=\frac{s}{4r^2}-\frac{r^2}{s} +\frac{1}{2s} +\frac{r^2}{4s^3}
+\left(-\frac{s}{8r^4}-\frac{3}{8r^2s}-\frac{3}{8s^3} -\frac{r^2}{8s^5}\right) L.
\]
Similarly, from $F_2=2fF_1+r$ one gets
\begin{align*}
 F_2 =\frac{s^2}{4r^3} &  - \frac{r^3}{s^2} -\frac{s^2}{16r^5} +\frac{1}{r}
 -\frac{1}{4r^3} +\frac{5r}{4s^2} -\frac{3}{8rs^2} +\frac{r^3}{2s^4}
 -\frac{r}{4s^4}-\frac{r^3}{16s^6} \\
 &+ LH_2 +L^2 J_2,
\end{align*}
with
\[
 -\frac{s^2}{r^5}< H_2 < 0 , \quad   -\frac{s^2}{r^9}< J_2 <0.
\]

The recurrence relation $F_{i+1} = 2 fF_i-F_{i-1}$ together with 
the chain of inequalities 
$s> r^{1.23} > 20^{0.23}r >31 r$ give polynomial expressions in $L$ of
the form
\begin{align*}
 F_3= & \frac{(16r^4 - 8r^2 + 1)s^3}{64r^8} + \frac{(32r^4 - 22r^2 + 3)s}{32r^6}
 + \frac{128r^4 - 96r^2 + 15}{64r^4s} \\
& {} +\frac{-16r^6 + 32r^4 - 26r^2 + 5}{16 r^2 s^3} +\frac{48r^4 - 56r^2 + 15}{64 s^5}
+ \frac{-6r^4 + 3r^2}{32 s^7} +\frac{r^4}{64 s^9}\\
& {} + L H_3+L^2 J_3 +L^3 K_3
\end{align*}
with
\[
-\frac{17s^3}{128r^6}< H_3 < 0 , \quad   -\frac{9s^3}{256r^{10}}< J_3 <0,
\quad -\frac{s^3}{256r^{14}}< K_3 <0,
\]
and 
\begin{align*}
 F_4= & \frac{(64r^6 - 48r^4 + 12r^2 - 1)s^4}{256r^{11}} + \frac{(32r^6 - 36r^4 + 11r^2 - 1)s^2}{32r^9}
\\
& {} + \frac{128r^6 - 192r^4 + 69r^2 - 7}{64r^7}  +\frac{96r^6 - 144r^4 + 60r^2 - 7}{32 r^5 s^2} \\
& {} +\frac{-128r^8 + 352r^6 - 504r^4 + 250r^2 - 35}{128 r^3 s^4}
+ \frac{32r^6 - 60r^4 + 39r^2 - 7}{32 rs^6} \\
& {} + \frac{-24r^5+ 27r^3 - 7r}{64 s^8} + \frac{2r^5 - r^3}{32 s^{10}} -\frac{r^5}{256s^{12}}\\
& {} + L H_4+L^2 J_4 +L^3 K_4  +L^4 M_4
\end{align*}
with
\[
-\frac{17s^4}{128r^7}< H_4 < 0 , \quad   -\frac{25s^4}{256r^{11}}< J_4 <0,
\]
\[
 -\frac{s^4}{128r^{15}} < K_4 <0, \quad   -\frac{s^4}{2048r^{19}}< M_4 <0.
\]

\medskip 

In view of Corollary~\ref{coFpo}, it is clear that Theorem~\ref{tede} 
is established as soon as we prove the next result.

\bepr{prsec}
Let $(1,b,c,d)$ be a  $D(-1)$--quadruple with $1<b<c<d$,
$b=r^2+1$, $c=s^2+1$, and $s=r^\theta$.
Then the following statements hold:

\emph{a)}  If $1.64 \le \theta <2.05$, then $|F_1| < 7\cdot 10^6$.

\emph{b)} If $1.40 \le \theta <1.64$, then $|F_2| < 3\cdot 10^6$.

\emph{c)} If $1.30 \le \theta <1.40$, then $|F_3| <  6\cdot 10^6$.

\emph{d)} If $1.23 \le \theta <1.30$, then $|F_4| <  2\cdot 10^6$.
\eepr 

In its proof we use an elementary fact, proved here for the sake of
completeness.

\bel{leteta}
Keep the notation from Proposition~\ref{prsec}. For $1.2< \theta <2.05$
one has
\[
 0 < \theta -\frac{\log c}{\log b} < \frac{1}{r^2}.
\]
\eel\bep 
The left inequality follows directly from Bernoulli's inequality.
Indeed, if $\eta =\log c/\log b $, then
\[
 c=r^{2\theta}+1 =(r^2+1)^\eta =r^{2\eta} \bigl( 1+r^{-2} \bigr)^\eta
 > r^{2\eta} \bigl( 1+ \eta r^{-2} \bigr) > r^{2\eta} +1.
\]

In view of the well-known fact $b^{-1} < \log b -\log r^2 <r^{-2}$,
the right inequality is consequence of $h(r^2)>\theta$, where
\[
 h(x)=\frac{1}{x}+\log x+\frac{x}{x^\theta+1}.
\]
The numerator of $h'$ is found to be 
$g(x)=(x-1)(x^\theta+1)^2+(1 -\theta) x^{\theta+2}+x^2$,
so that 
\begin{align*}
g''(x)=(4\theta^2+2\theta)x^{2\theta -1} & -(4\theta^2-2\theta)x^{2\theta -2} 
-(\theta^3+2\theta^2 -\theta -2)x^\theta \\
& +(2\theta^2+4\theta)x^{\theta -1}-(2\theta^2-2\theta)x^{\theta -2} +2.
\end{align*}
The sum of the last three terms in the above expression is obviously
positive and it is easily checked that the same is true for the sum of 
the other terms. Therefore, for $x>2$ and $\theta <2.05$ one has  
$g'(x)> g'(2)> g'(1) =13-\theta -\theta^2>0$ and  
$g(x) >g(2)>(6-3\theta)2^\theta +\theta +6>0$. We conclude that the 
function $h$ is increasing, so $h(r^2) > 13\log 10 >\theta$.
\eep 

\noindent
\emph{Proof of Proposition~\ref{prsec}.}
a)
As $\theta <2.05$, we can bound from above $F_1$ as follows:
\[
 F_1 < \frac{s}{4r^2}-\frac{r^2}{s} +\frac{1}{2s} +\frac{r^2}{4s^3}
 < \frac{s}{4r^2} <0.25 \, r^{0.05}.
\]
By Theorem~\ref{prale}, for $1.64\le \theta < 2.05$ it holds 
$b< 10^{50}$. Therefore, 
\[
 F_1 < 0.25 \, \bigl( 10^{25} \bigr)^{0.05} < 5.
\]

We bound from below $F_1$  quite similarly:
\[
 F_1> \frac{s}{4r^2}-\frac{r^2}{s} +\frac{1}{2s} +\frac{r^2}{4s^3}
-\frac{s}{8r^4}-\frac{3}{8r^2s}-\frac{3}{8s^3} -\frac{r^2}{8s^5}
> -\frac{r^2}{s} \ge -r^{0.36}.
\]
Our program for computation of an absolute upper bound for $b$ 
iterates over $\log c/\log b$ not over $\theta$, which explains the 
need for Lemma~\ref{leteta}. The computations show $b<10^{38}$ 
when $\theta$ is in the range $1.639< \theta <2.05$, so that
\[
 F_1> -\bigl( 10^{19} \bigr)^{0.36} > -7\cdot 10^6.
\]

b) Under the current hypothesis we get 
\[
F_2 < \frac{s^2}{4r^3}   - \frac{r^3}{s^2}  +\frac{2}{r} < \frac{s^2}{4r^3} 
= 0.25 \, r^{2\theta -3}< 0.25 \, r^{0.28}.
\]
Using the bound on $b$ stated in Theorem~\ref{prale} results
in a bound on $F_2$ outside the desired range. Therefore we 
split the interval where $\theta$ takes its values.
From the output of our program for computation of an absolute upper 
bound for $b$ we see that  $b<10^{49.3}$ when  $1.499< \theta <1.64$, 
so that 
\[
 F_2 < 0.25 \, \bigl( 10^{24.65} \bigr)^{0.28} < 3\cdot 10^6 \quad \mbox{when {$1.5 \le  \theta <1.64$}}.
\]
For $1.40\le \theta <1.50$ one  gets at once
\[
 F_2< 0.25 \, r^0<1.
\]

We similarly  bound $F_2$ from below:
\[
 F_2 >\frac{s^2}{4r^3}   - \frac{r^3}{s^2} -\frac{s^2}{16r^5}
 -\frac{s^2}{r^5} -\frac{s^2}{r^9} >  - \frac{r^3}{s^2} = -r^{3-2\theta}.
\] 
When $1.5\le \theta <1.64$, this gives 
\[
F_2 > -1,
\]
while on the subinterval $1.40\le \theta <1.50$ it implies
\[
 F_2 > -r^{0.2} > -\bigl( 10^{31.25} \bigr)^{0.2} > -2\cdot 10^6
\]
because $b<10^{62.5}$ on this subinterval.

c)  From the expression for $F_3$ we first obtain
\[
 F_3< \frac{s^3}{4r^4}+\frac{s}{r^2} +\frac{2}{s}-\frac{r^4}{s^3} +\frac{2r^2}{s^3}
 +\frac{3r^4}{4s^5}  <\frac{s^3}{4r^4}+\frac{2s}{r^2}-\frac{r^4}{s^3}<\frac{s^3}{4r^4}.
\]
As seen from Theorem~\ref{prale}, one has $b< 6.26\cdot 10^{73}$ when 
$ \theta \ge 1.299$. Therefore, the upper bound for $F_3$ just obtained 
can be bounded from above as follows: 
\[
 0.25 \, r^{0.2}< 0.25 \, \left(62.6^{0.5}\cdot 10^{36}  \right)^{0.2} < 6\cdot 10^6.
\]

To obtain a lower bound for $F_3$, we can ignore all fractions but the first, the fourth
and the sixth in its free term and replace the coefficients of positive powers of
$L$ by their respective lower bounds.  We thus get
\[
F_3 > \frac{3s^3}{16r^4}+\frac{s}{2r^2}-\frac{r^4}{s^3}-\frac{3r^4}{16s^7}
-\frac{17s^3}{128r^6} -\frac{9s^3}{256r^{10}} -\frac{s^3}{256r^{14}} 
> \frac{s}{2r^2}-\frac{r^4}{s^3}-\frac{3r^4}{16s^7} > -\frac{r^4}{s^3},
\]
whence
\[
 F_3> -r^{0.1} > -\left(62.6^{0.5}\cdot 10^{36}  \right)^{0.1} > -5000.
\]

d) We similarly  see that it holds
\[
 F_4<  \frac{s^4}{4 r^5} < 0.25 \,  r^{0.2} 
 < 0.25 \,  \left(10^{34.5} \right)^{0.2} < 2 \cdot 10^6
\]
and 
\[
F_4>  -\frac{r^5}{s^4} > -r^{0.08} > - \left(10^{34.5} \right)^{0.08} > -600.
\]

Proposition~\ref{prsec} being established, the proof of the nonexistence
of $D(-1)$--quadruples is complete.


\medskip

\textbf{Acknowledgments.} The results reported in this paper 
would not have been obtained without the computations performed 
on the computer network of IRMA and Department of Mathematics
and Computer Sciences of Universit\'e de Strasbourg. The authors
are grateful to Ryotaro Okazaki and Yasutsugu Fujita
for drawing attention to Matveev's ignored result from~ \cite{Mat98}.

\end{document}